\documentclass[11pt]{article}
\usepackage{amsmath} 
\usepackage[e]{esvect}
\usepackage{graphicx}
\usepackage[adobe-utopia]{mathdesign}
\usepackage[T1]{fontenc}
\usepackage[english]{babel}
\usepackage{diagbox,booktabs}
\usepackage{color}

\usepackage{hyperref}

\title{Nearly neighbourly families of standard boxes}
\date{}
\author{Jacek Bojarski, Andrzej P. Kisielewicz
and Krzysztof Przes{\l}awski}
\newtheorem{pr}{\sc Proposition}
\newtheorem{lemat}[pr]{\sc Lemma}
\newtheorem{tw}[pr]{\sc Theorem}
\newtheorem{wn}[pr]{\sc Corollary}
\newtheorem{df}{\sc Definition}

\newtheorem{uw}{\sc Remark}
\newtheorem{uwi}[uw]{\sc Remarks}
\newenvironment{uwa}{\begin{uw} \rm}{\end{uw}}

\newtheorem{nap}{\sc Example}
\newenvironment{npr}{\begin{nap} \rm}{\end{nap}}
\newtheorem{nps}[nap]{\sc Examples}

\def\ka #1{\mathscr{#1}}
\def\kal #1 #2{\mathscr{#1}^{#2}}
\def\proof{\noindent \textit{Proof.\,\,\,}}
\def\skok{\vskip 0.1in\noindent}
\def\Skok{\vskip 0.15in}
\def\enn{\mathbb{N}}
\def\zet{\mathbb{Z}}

\def\er{\mathbb{R}}

\def\aut #1{\operatorname{Aut} (#1)}
\def\Aut #1 #2{\operatorname{Aut}^{#1} (#2)}

\def\supp #1{\operatorname{supp}\, #1}

\def\bred #1 {\colorbox{red}{ #1}}
\def\red #1 {{\color{red} #1 }} 
\def\dkrop{\colon\hspace{-2pt} }

\graphicspath{{figpdfkolor/}}
\DeclareGraphicsExtensions{.pdf,.png} 

\begin{document}
\maketitle
\begin{abstract}
It is demonstrated that each nearly neighbourly family of standard boxes in $\er^3$ has at most 12 elements. 
A combinatorial classification of all such families that have exactly 12 elements
is given.  All families satisfying an extra property called incompressibility are described. Compressible families are discussed briefly.

\end{abstract}

\section{Introduction}

A  \textit{standard $n$-box}  or an $n$-\textit{interval}  is the Cartesian product of $n$ ordinary closed intervals of positive length. Two $n$-intervals $I=I_1\times\cdots\times I_n$ and $J=J_1\times\cdots\times J_n$ are adjacent if there is $i \in [n]$ such that  $I_i$ and $J_i$ have exactly one point in common. Both the family and  the infinite graph of all $n$-intervals with  the adjacency just defined are denoted by $\ka I^n$.  This convention extends to any subfamily  $\ka G$ of $\ka I^n$: The  same symbol $\ka G$ is for the graph with adjacency inherited from $\ka I^n$. A subfamily $\ka G\subset \ka I^n$ is \textit{nearly neighbourly} if  it is a clique in $\ka I^n$.  

The main purpose of the present investigation is to demonstrate that the maximum cardinality of a nearly neighbourly family in $\ka I^3$ is 12 (Theorem \ref{be3}), and to give a more or less full description of all such families from both a combinatorial and a geometric points of view.  It should be pointed out that the case of the so-called compressible cliques (see Section 6) is addressed rather superficially. Only results without proofs are  presented.  A more detailed analysis will possibly be published elsewhere.

The notion of  a nearly neighbourly family of intervals is a specialization of  a more general concept promoted by  Zaks \cite{Zaks1}: A family  $\ka P$ of  $n$-dimensional  convex polytopes in $\er ^n$ is said to be \textit{nearly neighbourly}, if for every two polytopes $P$, $Q$ belonging to $\ka P$ there is a hyperplane separating them  that contains a facet of $P$  and a facet of  $Q$.  In fact, researchers have paid more attention  to neighbourly families. Let us remind that $\ka P$ is \textit{neighbourly} if the intersection of any two members of $\ka P$   is of dimension $n-1$. In general, such a $\ka P$ can be of arbitrary finite cardinality unless {$n\le 2$};  in addition, one may even assume that the members of $\ka P$ are affinely equivalent \cite{Ku} or congruent \cite{EK}. However,  if a nearly neighbourly  family $\ka P$ consists of polytopes having their number of facets bounded from the above by $m$, then, as mentioned in \cite{Zaks2}, 
\begin{equation}
\label{Perlest}
 |\ka P|\le 2^m.  
\end{equation}
The proof is based on an idea of Perles  \cite{Perles}.  (A similar method has been employed by other researches even earlier; see \cite{Furedi} for further details).  It is conjectured that if $\ka P$ consists of $n$-dimensional simplices, then the above estimate can be improved by factor 2; that is, $|\ka P|\le2^n$.  The conjecture is open for all $n\ge 3$.  (In case of $n=3$, the best known estimate is 14 (see \cite{Furino})). It is even open for neighbourly families of simplices in dimension 4.  As it concerns tetrahedra, it was  Bagemihl \cite{Bagemihl} who raised the question. He constructed a neighbourly family of 8 tetrahedra, and  speculated that it is a family of maximum cardinality. Subsequently, Baston \cite{Baston} proved that a neighbourly family of tetrahedra has at most 9 elements. And finally  Zaks \cite{Zaks3} was able to show that it has at most 8 elements as expected.  His work depends heavily on Baston's research and the Graham---Pollak theorem \cite{GP}  on minimal biclique partitions of complete graphs.  It was also Zaks \cite{Zaks4} who constructed a neighbourly family of $n$-dimensional simplices consisting of $2^n$ members for $n> 3$.  

We know  just a couple of papers devoted to (nearly) neighbourly families of standard boxes. Zaks \cite{Zaks5}  proved that the maximum cardinality of a neighbourly family of $n$-intervals is $n+1$. Again, the proof  depends on the Graham---Pollak theorem.  In \cite{Alon}, N. Alon studied  $k$-neighbouring families of $n$-intervals. Let us remind that $\ka P$ is such a family if for every two members $P$ and $Q$  of $\ka P$ one has $n-k\le \dim P\cap Q \le d-1$.  He gave estimates from below and above for the maximum cardinality of a $k$ neighbourly family of $n$-intervals. There are two works  \cite{Simon1, Simon2} by J. D. Simon on (nearly) neighbourly families of quadrilaterals. Some of her results  will be discussed in Section \ref{maxc 2D}.    

We  begin with showing that for families of $n$-intervals  Perles' estimate (\ref{Perlest}) can be slightly improved (Proposition \ref{rekurencja}, Remark \ref{est}).  In this context, we introduce a bulk of notions instrumental for further presentation. 

Let $X$ be  a non-empty set. We denote by $\enn X$ the abelian semigroup of all finite formal sums of elements of $X$.  (The elements of $\enn X$ will also be called   \textit{combinations}). Every combination $\gamma \in \enn X$  is uniquely determined by a function $k\colon X\to \enn$ with support $\{x\colon k_x>0\}$ of finite cardinality. We shall use the following notation related to $\gamma$: 
$$
\gamma=\sum_{x\in X} k_xx=\sum k_x x, \quad |\gamma|=\sum {k_x}, \quad  \supp \gamma =\{x\colon k_x>0\}.
$$
On several occasions we will use a naturally defined inequality relation in $\enn X$: For  $\delta=\sum l_xx$ and $\gamma=\sum k_x x$,  we write $\delta\le \gamma$  if and only if  $l_x\le k_x$, for every $x\in X$.  If $\delta\le \gamma$,  then we say that $\delta$ is  a \textit{subcombination} of $\gamma$.  If $l_x\in \{k_x,0\}$, for every $x\in X$, then   $\delta$ is  an \textit{induced} subcombination of $\gamma$. 
Suppose  $G$  is a graph whose set of vertices $V(G)$ equals $X$. Then each $\gamma$ determines in a natural way the graph $G_\gamma$:
\begin{eqnarray} 
V(G_\gamma)&=& \bigcup _{x\colon k_x>0} \{ (x,1),\ldots, (x,k_x)\}; \nonumber\\
E(G_\gamma)&=&\{ \{(x,i),(y,j)\}\colon \{x,y\}\in E(G),  1\le i\le k_x, 1\le j\le k_y\}.\nonumber
\end{eqnarray}
In what follows,  we shall often appeal to $\gamma$ itself as to $G_\gamma$.  Consequently, we shall write $V(\gamma)$, $E(\gamma)$ rather than $V(G_\gamma)$, $E(G_\gamma)$. 
\begin{pr}
\label{lancuchy}
Let $G$ be a graph and let  $\gamma =\sum k_x x$ be an element of $\enn V(G)$. Then
$$
|V(\gamma) |=\sum k_x , \qquad  |E(\gamma)|=\sum_{\{x,y\}\in E(G)} k_xk_y.
$$
\end{pr}

The \textit{clique number} $\omega(G)$ of a graph $G$ is the cardinality of a maximum clique contained in $V(G)$. If $\gamma\in \enn V(G)$,  Then
$\omega(\gamma)=\omega(G_\gamma)$ is equal to $\omega (G[\supp \gamma])$, the clique number of the subgraph of $G$ induced by $\supp \gamma$. The  \textit{independence number}  $\alpha(G)$ of $G$ is the clique number of the complement graph  of $G$.  

Suppose that two non-empty sets $X$ and $Y$, and a mapping $f\colon X\to Y$ are given. Then $f$ induces the mapping $f_*\colon\enn X\to \enn Y$ 
$$
f_*(\gamma) =\sum k_x f(x).
$$ 
\begin{pr}
Given two graphs $G$ and $H$, and a homomorphism of graphs $f\colon V(G)\to V(H)$.  Then for every $\gamma\in \enn V(G)$, the graph $f_*(\gamma)$ is a homomorphic image of  $\gamma$.  In particular,  $\alpha(f_*(\gamma))\le \alpha(\gamma)$.
\end{pr} 

Our next proposition is rather obvious.
\begin{pr}
\label{or}
Let $G_1$ and  $G_2$ be two graphs and let $G= G_1 * G_2$ be their \textit{disjunctive product}; that is, $V(G)=V(G_1)\times V(G_2)$ and $\{(u_1,u_2),(v_1,v_2)\}\in E(G)$ if and only if $\{u_1,v_1\}\in E(G_1)$ or $\{u_2,v_2\}\in E(G_2)$. Let  $C$ be a finite clique in $G$ and  $\gamma_i=\sum_{x\in C} x_i$, $i=1,2$. Then 
$$
\alpha(\gamma_1)\le \omega(\gamma_2).
$$
  \end{pr} 
\proof
Let $D\subseteq C$ be a set of maximum cardinality such that $\sum_{x\in D} x_1$ is independent. Then  $\sum_{x\in D}x_2$ is a clique and $|D|=\alpha(\gamma_1)$. Therefore,
$$
 \alpha(\gamma_1)= \left|\sum_{x\in D}x_2\right| \le\omega(\gamma_2) .
$$
\hfill$\square$   
\begin{uwa}
\label{izo2}
Let $K=\{i_1<i_2<\ldots <i_k\}$ be a proper subset of $[n]$. Let $K^\mathsf{c}$ be its complement.  For $I \in \ka I^n$, let $I_K=I_{i_1}\times\cdots \times I_{i_k}$. The mapping $I\mapsto (I_{K}, I_{K^\mathsf{c}})$ defines an isomorphism  between graphs $\ka I^n$ and $\ka I^{|K|}*\ka I^{|K^\mathsf{c}|}$.
\end{uwa}
By this remark and the preceding proposition we have
\begin{wn}
\label{wnor}
Let $K$ be a proper subset of $[n]$. Let $\ka C$ be a clique in $\ka I^n$. Let $\gamma_L=\sum_{I \in \ka C} I_L$, where $L=K,K^\mathsf{c}$. Then
$$
\alpha(\gamma_{K})\le \omega(\gamma_{K^\mathsf{c}}).
$$ 
\end{wn} 

In what follows, the notation of Corollary \ref{wnor} is  employed in a more general setting: Let $\gamma=\sum_I k_I I\in \enn\ka I^n$, and $K$ be a proper subset of $[n]$; then $\gamma_K =\sum_I k_I I_K$.   If $K=\{i\}$, then we shall often write $\gamma_i$  rather than $\gamma_{\{i\}}$. 
  	
\section{Incompressibility. An upper bound for the cardinality of a maximum clique in \texorpdfstring{$\ka I^n$}{dimension n}\label{upper bound}} 

For $s \in \enn$, let $\ka I(s)=\{[0,1], [s,s+1]\}\cup \{[i,j]\colon 1\le i<j\le s\} $. 
Let $\ka G$ be a finite subfamily of $\ka I^1$.  By applying an appropriate homeomorphism  $H$ of $\er$, one can transform $\ka G$  onto $\ka G'=\{ I'\colon I'=H(I),\,  I \in  \ka G\}$  so that  the endpoints of the intervals $I'\in \ka G'$ are positive integers.  Therefore, there is an $s$ and a homomorphism of graphs
$f\colon \ka G\to \ka I(s)$.  
{We shall be concerned with properties of the family $f(\ka G)$ for minimal $s$.}    

\begin{pr}
\label{shom}
 Let a finite nonempty family $\ka G\subseteq \ka I^1$ be given. Let $s=s(\ka G)$ be the minimum number for which  there is a  homomorphism of graphs $f\colon \ka G \to\ka I(s)$. Then 
\begin{itemize}
\item[$(1)$] 
$\{[i,i+1]\colon 0\le i\le s\}\subseteq f(\ka G);$
\item[$(2)$]
$\alpha(\ka G)\ge \left\lfloor\frac s 2\right\rfloor+1;$
\item[$(3)$]
$\{[1,3], [s-2,s]\}\subseteq f(\ka G)$, whenever $s\ge 3.$
\end{itemize}
\end{pr}
\proof
 There is nothing to prove if $s=0$. The case $s=1$ is rather obvious. Moreover, $s\not=2$  as the mapping $[2,3]\mapsto [0,1]$, $[i,i+1]\mapsto [i,i+1]$ for $i=0, 1$ is a homomorphism of $\ka I(2)$ onto $\ka I(1)$. Therefore we may assume $s\ge 3$.  Suppose $[0,1]\not\in f(\ka G)$.  For every $I=[a,b]$ belonging  to $f(\ka G)$, let as set 
$$
I'=
\left\{
\begin{array}{ll}
\text{$[1, b-1]$}    &  \text{if $a=0$ or $a=1$, $b>2$,}\\
\text{$[a-1,b-1]$} &  \text{if $a=1$,  $b=2$ or $a>1$.}\\
\end{array} 
\right.
$$
The mapping  $g$ defined by $I'=g(I)$ is a homomorphism of $f(\ka G)$ into $\ka I(s-1)$. Therefore, the composite  $g\circ f$ maps $\ka G$ into $\ka I(s-1)$ contrary to the definition of $s$.  
Since the mapping $[a,b]\stackrel h\mapsto [s+1-b, s+1-a]$ is an automorphism of $\ka I(s)$ which sends $[s,s+1]$ on $[0,1]$,   the interval $[s,s+1]$ has to belong to $f(\ka G)$ also.

Suppose now that  $[i,i+1]\not \in  f(\ka G)$, where   $i\in [1,s-1]$.  Then define the correspondence $I=[a,b]\stackrel g\mapsto I'$ as follows
$$
I'=
\left\{
\begin{array}{ll}
\text{$[a, b]$}    &  \text{if $b\le i$,}\\ 
\text{$[a,b-1]$} &  \text{if $a\le i$,  $b\ge i+1$,}\\
\text{$[a-1,b-1]$} &  \text{if $a\ge i+1$.}\\
\end{array} 
\right.
$$
Again, the composite $g\circ f$ is a homomorphism  of $\ka G$ into  $\ka I(s-1)$ contradicting the minimality of $s$. 

In order to prove $(2)$, it suffices to observe that by $(1)$, the independent set of intervals $\{[2i,2i+1]\colon 0\le i\le \left\lfloor  \frac s 2\right\rfloor \}$ is contained in $f(\ka G)$,  and that $\alpha(f(\ka G)) \le \alpha(\ka G)$. 

As  it concerns  the intervals $[1,3]$ and $[s-2,s]$, the  automorphism $h$ transposes them.  Moreover, there is an automorphism (see Appendix A)  which transposes $[1,3]$ with $[2,3]$. Therefore, both $[1,3]$ and $[s-2,s]$ have to belong to $f(\ka G)$.  \hfill $\square$
\begin{uwa}
Suppose $\ka G\subset \ka I^1$ consists of intervals of odd length with integer endpoints.  Then  the mapping $[a,b]\mapsto [a\operatorname{mod\, 2}, 1+ (a\operatorname{mod\, 2)} ]$ is a graph homomorphism from $\ka G$ to $\ka I(1)$.   If there are two intervals in $\ka G$ with right endpoints of different parity, then $s(\ka G)=2$ and consequently $\ka G$ is bipartite. If all intervals in $\ka G$ have their right endpoints of equal parity, then $s(\ka G)=0$ and $\ka G$ has no edges. 
\end{uwa}

A subgraph $\ka G\subseteq\ka I^1$ is \textit{incompressible} if there is $s$ such that  $\ka G\subseteq \ka I(s)$ and  there is no a homomorphism  $f\colon \ka G \to \ka I(s')$, where  $s'<s$.  Otherwise, $\ka G$ is called compressible.   A combination $\gamma \in   \enn \ka I(s)$  is incompressible  (compressible) if $\supp \gamma$ is an incompressible (compressible) subgraph of $\ka I(s)$.   As an immediate consequence of the preceding proposition we have
\begin{pr}
\label{irred}
If $\gamma \in   \enn \ka I(s)\setminus \{0\}$ is incompressible  then  
\begin{itemize}
\item[$(1)$]
$s \le 2\alpha(\gamma)-1;$  
\item[$(2)$] 
$\ka S(s)=\{[i,i+1]\colon  0\le i\le s \}\cup \{ [1,3], [s-2, s]\}$ is contained in $\supp\gamma$, whenever $s\ge 3$.
\end{itemize}
\end{pr}

\begin{lemat}
\label{inne}
Let $\gamma\in \enn\ka I^1$ be given. If $\alpha(\gamma)>1$, then
$$
|\gamma|\le 4\alpha(\gamma)-3.
$$
\end{lemat}
\proof
We may assume that $\supp \gamma \subseteq \ka I(s)$ for some $s$. 
If $\gamma$ is compressible, then let us choose a new $s$ and  a homomorphism $f\colon \supp \gamma \to \ka I(s)$ so that $\gamma'=f_*(\gamma)$ is incompressible. 
If $s=0$ or $s=1$, then since $\alpha(\gamma)>1$, we obtain
$$
|\gamma|=|\gamma'|\le 2\alpha(\gamma') \le 2\alpha(\gamma)\le 4\alpha(\gamma)-3.
$$
As $\ka I(2)$ is compressible, we may assume that $s>2$.  Then, by the preceding lemma, $[0,1]$ and $[2,3]$ are elements of $\supp \gamma'$ and consequently $\alpha(\gamma')> 1$.  Henceforth, we may further assume that  $\gamma\in \enn\ka I(s)$ is incompressible  and $s>2$.  Let us set $B=\{0,1\}^s$ and 
$$
B(I)=
\left\{
\begin{array}{ll}
\{1\}\times \{0,1\}^{s-1}, & \text{if $I=[0,1]$}\\
\{0,1\}^{i-1}\times \{0\}\times \{0,1\}^{j-i-1}\times \{1\}\times \{0,1\}^{s-j},  & \text{if  $I=[i,j]$ and $1\le  i <j\le s$}\\
  \{0,1\}^{s-1}\times\{0\}, & \text{if $I=[s,s+1]$}\\ 
\end{array}
\right.
$$
Let  $\gamma^\sharp\colon B\to \er$ be associated with $\gamma=\sum_{I\in \ka I(s)} k_I I$ by the formula
$$
\gamma^\sharp(x)=\sum_{I\in \ka I(s)} k_I \mathbf{1}_{B(I)}(x). 
$$
Observe that  $\gamma^\sharp$ is bounded from the above by  $\alpha(\gamma)\mathbf{1}_B$.  Moreover, these two functions do not coincide, as if   $x_0=(0,\ldots,0)$,  then,  by  Proposition \ref{shom},
$k_{[0,1]}>0$ and consequently
$$
 \gamma^\sharp(x_0)=k_{[s,s+1]}<k_{[0,1]}+k_{[s,s+1]}\le \alpha(\gamma).
$$
Therefore, summing  up these functions with respect to $x$ yields
$$
\alpha(\gamma) |B|> \sum_I k_I |B(I)|.
 $$
If $I=[0,1]$ or $I=[s,s+1]$, then  $|B(I)|=\frac 1 2 |B|$;  if not, $|B(I)|=\frac 1 4  |B|$.  By Proposition \ref{shom},  $k_{[0,1]}$ and $k_{[s,s+1]}$ are both greater or equal to $1$, thus 
$$
\alpha(\gamma)> \frac{k_{[0,1]} + k_{[s,s+1]}} 4 +\frac 1 4 \sum_{I\in \ka I(s)} k_I\ge \frac 1 2 +\frac{|\gamma|}4.
$$
\hfill $\square$

We denote by $b_m$ the maximum cardinality of a clique in $\ka I^m$. 

\begin{pr}
\label{oszac}
  Let $K$ be a proper subset of $[n]$. Let $\ka C$ be a clique in $\ka I^n$. Let $\gamma_{K}=\sum_{I \in \ka C} I_{K}$. Then
$$
\alpha(\gamma_{K})\le b_{n-|K|}.
$$ 
\end{pr}
\proof
By Corollary \ref{wnor} and the fact that $\supp (\gamma_{K^\mathsf{c}})\subset\ka I^{n-|K|}$,  one has
$$
\alpha(\gamma_K)\le \omega(\gamma_{K^\mathsf{c}})=\omega(\supp(\gamma_{K^\mathsf{c}}))\le \omega\left(\ka I^{n-|K|}\right)=b_{n-|K|}.
$$
\hfill $\square$ 
\begin{pr} 
\label{rekurencja}
For every $n\ge 2$
$$
b_n\le 4b_{n-1}-3.
$$
Moreover, $b_2=5$. 
\end{pr}

\proof
Suppose $\ka C$ to be a finite clique in $\ka I^n$.  Let $\gamma_1=\sum_{I\in \ka C} I_1$.
By Lemma \ref{inne} and Proposition~\ref{oszac}, one gets
$$
|\ka C|=|\gamma_{1}|\le 4\alpha(\gamma_{1})-3\le 4b_{n-1}-3,
$$
whenever $\alpha(\gamma_{1})>1$.  If $\alpha(\gamma_{1})=1$ , then $|\gamma_{1}|\le 2$  and the inequality $|\ka C|\le 4b_{n-1}-3$ is still valid. 

To prove the second part, observe that since $b_1=2$, it follows by the first part that $b_2\le 5$. Any of the configurations depicted in Figure \ref{piatki} shows that $b_2$ is at least $5$. \hfill $\square$ 

\begin{uwa}
The fact that $b_2=5$ is due to Simon \cite[Theorem 5.7]{Simon1}.  Her proof seems to be different.  It depends on the analysis of the so-called Baston matrix of a maximum family of 2-intervals.
\end{uwa}   
\begin{wn}
Let $k\ge 1$ and $n\ge k$. Then  
$$
b_n\le \frac{b_k-1}{4^k}4^{n}+1.
$$
\end{wn}

\begin{uwa}
\label{est}
As we have announced, we shall demonstrate later on that $b_3=12$. Therefore, for every $n\ge 3$
$$
b_n\le \frac{11}{64}4^{n}+1.      
$$
\end{uwa}
\begin{uwa}
\label{superm}
If $\ka C$ is a clique in $\ka I^k$ and $\ka D$ is a clique in $\ka I^l$, then $\ka E=\{I\times J \colon I \in \ka  C, J\in \ka D\}$ is  a clique in $\ka I^{k+l}$. Therefore, we have
$$
b_kb_l\le b_{k+l}.
$$ 
In particular, for every $n>1$ 
$$
b_n\ge 2b_{n-1},
$$
as $b_1=2$.
\end{uwa}
\section{Maximum cliques in \texorpdfstring{$\ka I^2$}{dimension 2} \label{maxc 2D}}

We are going to describe all cliques of maximal cardinality in $\ka I^2$.  Let $\ka C$ be such a clique. By Proposition \ref{rekurencja}, $|\ka C|=5$. Let   
$E_{i}=\{\{I,J\}\in E(\ka C) \colon \text{$I_i$ and $J_i$ are adjacent in $\ka I^1$ }\}$, $i=1, 2$.  Let the graphs $G_i$  be defined by the equations $V(G_i)=\ka C$, $E(G_i)=E_i $.  Let  $\gamma_i =\sum_{I\in\ka C} I_i$.  It is clear that $G_i$ is isomorphic to the graph associated with $\gamma_i$.  Since $\ka I^1$ does not contain triangles, $\gamma_i$ and consequently $G_i$ cannot contain triangles as well.  By the definition of adjacency in  $\ka I^2$, the sets $E_1, E_2$ cover $E(\ka C)$. These facts imply that the graphs $G_i$ have no vertices of degree higher than two and $|E_i|\le 5$.   As $\ka C$ is a 5-clique, its set of edges $E(\ka C)$ has exactly ten elements. Therefore,  $\{E_1, E_2\}$ is a partition of $E(\ka C)$ and $|E_i|=5$ for $i=1,2$. It is rather clear now that $G_1, G_2$ gives us a decomposition of $\ka C$ into two 5-cycles.   

Given a strictly increasing sequence of numbers $a_0,\ldots, a_4$. Then
$$
\Gamma = \{[a_k,a_{k+1}]\colon k=0,\ldots , 3\}\cup \{[a_1,a_3]\}.
$$
is a 5-cycle in $\ka I^1$. It is easy to verify, that all 5-cycles in $\ka I^1$ are of this form.  Let  $\gamma\in \enn\ka I^1$, be a 5-cycle. It is equally easy to  observe that there is a unique 5-cycle $\Gamma\in \ka I^1$ such that $\gamma=\sum_{I\in \Gamma} I$.  Henceforth, we  have just reduced the problem of describing all cliques of maximal cardinality in $\ka I^2$ to the following \textit{  construction} problem:
\begin{quotation}
Given two $5$-cycles $\Gamma_i$, $i=1,2$ in $\ka I^1$. Find all possible  5-cliques $\ka C\subset \ka I^2$ such that 
\begin{equation}
\label{five}
\ka C_i:=\{I_i \colon I \in \ka C\}=\Gamma_i.
\end{equation}
\end{quotation}  

\begin{figure}
\centering
\includegraphics{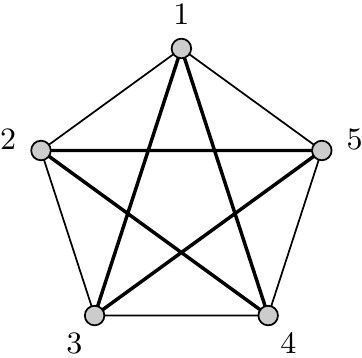}
\caption{$G_1=12345;\, G_2=13524$.\protect\label{5klika}}
\end{figure}

Let $a_{i,0},\ldots, a_{i,4}$,  be the sequence defining the 5-cycle $\Gamma_i$ for $i=1,2$.  We can  label each interval belonging to $\Gamma_1$ as follows:
$$
\text{$[a_{1,k},a_{1,k+1}]\mapsto k+1$ for $k=0,\ldots, 3$ and $[a_{1,1},a_{1,3}]\mapsto 5$.}  
$$
Since $I\mapsto I_1$ is a one-to-one mapping from $\ka C$ onto $\Gamma_i$ for every $\ka C$ to be constructed, it determines a labelling of $\ka C$ by composition. Now, we can identify $\ka C$ with the complete graph on $\{1,2,3,4,5\}$ (see Fig. \ref{5klika}), where the edges of the cycle $G_1$ corresponding to $\Gamma_1$ are the  edges of the pentagon   while the edges of the other cycle $G_2$ corresponding to $\Gamma_2$ are the edges of the pentagram.  

\begin{figure}[hbp]
\centering
\def\wsp{0.28}
\begin{tabular}{ccc}
\includegraphics[width=\wsp\textwidth]{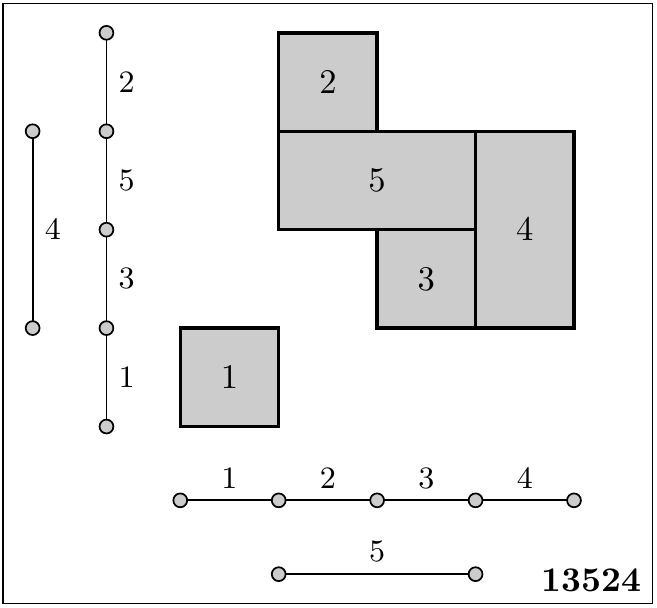} &
\includegraphics[width=\wsp\textwidth]{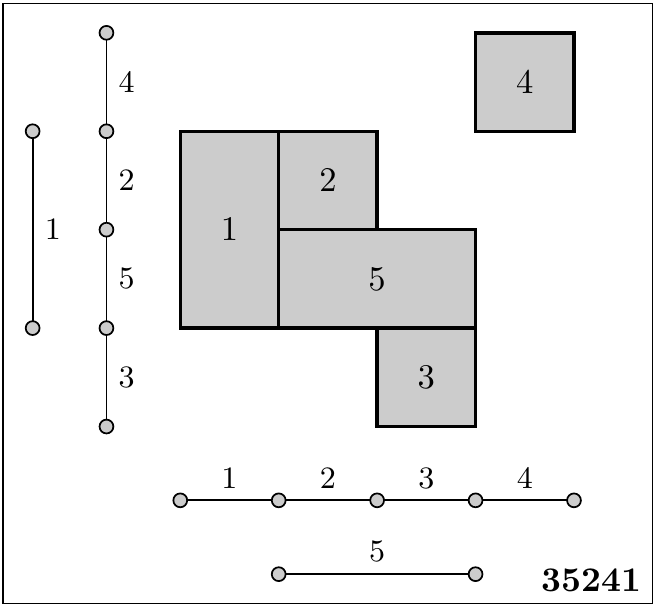}&
\includegraphics[width=\wsp\textwidth]{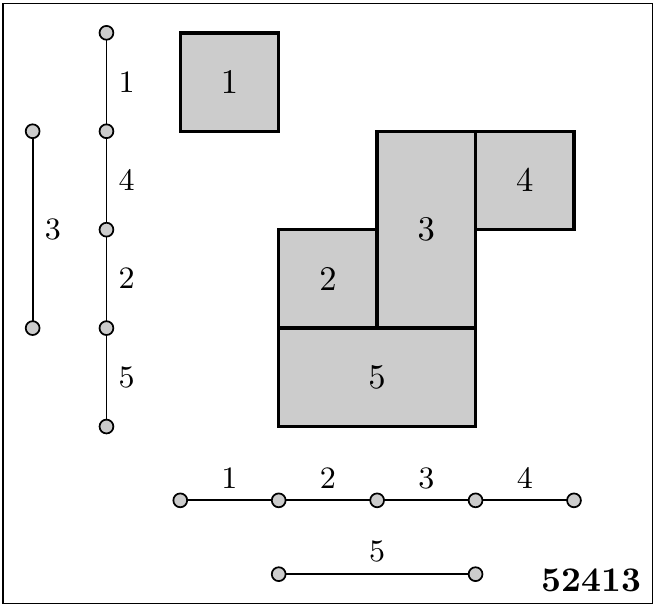}\\ 
\includegraphics[width=\wsp\textwidth]{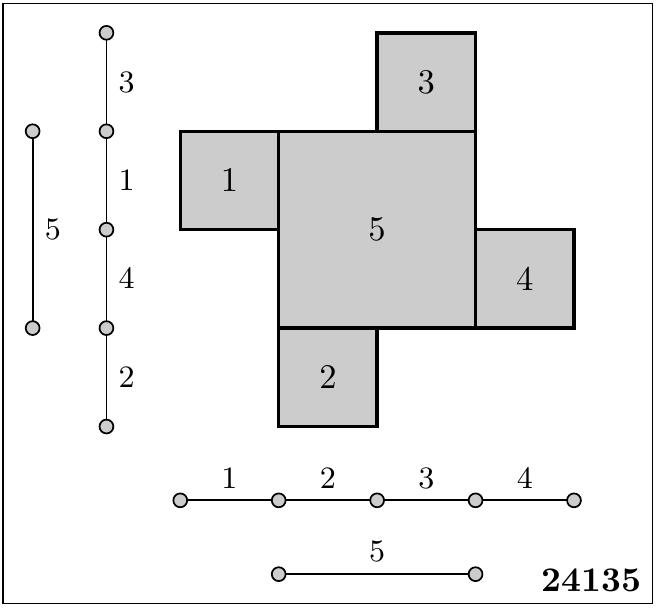} &
\includegraphics[width=\wsp\textwidth]{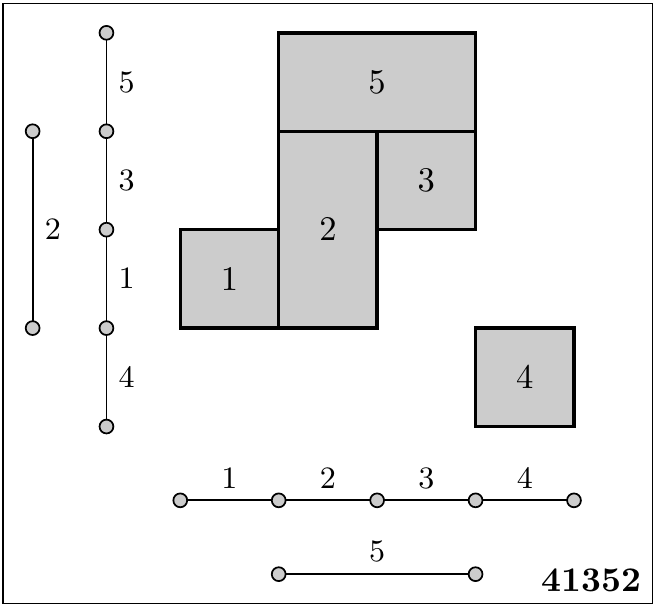}&
\includegraphics[width=\wsp\textwidth]{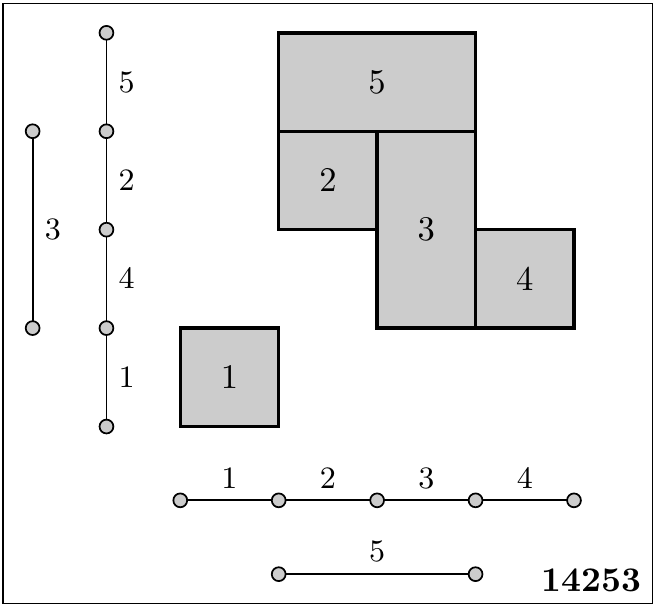}\\ 
\includegraphics[width=\wsp\textwidth]{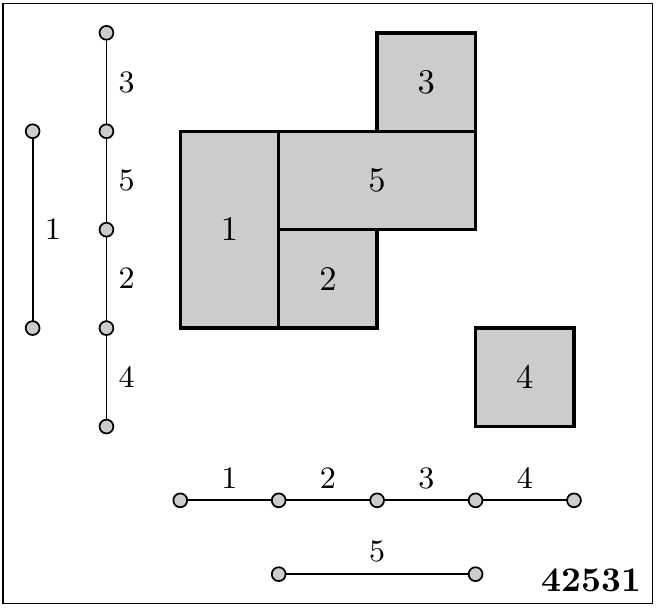} &
\includegraphics[width=\wsp\textwidth]{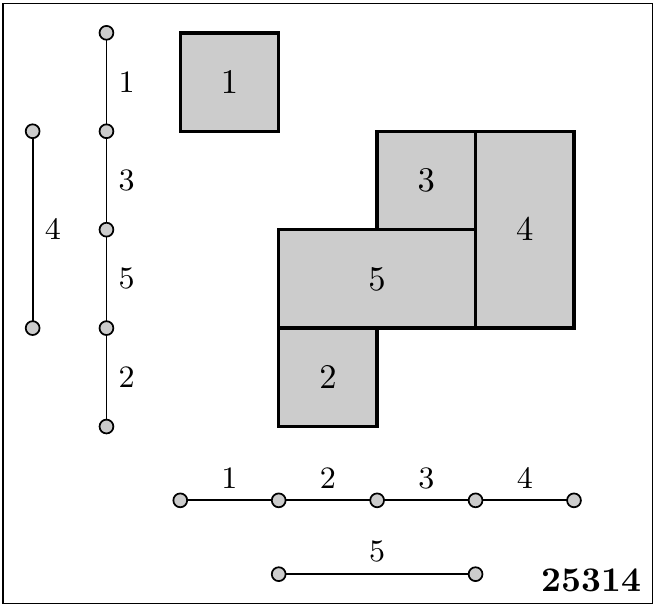}&
\includegraphics[width=\wsp\textwidth]{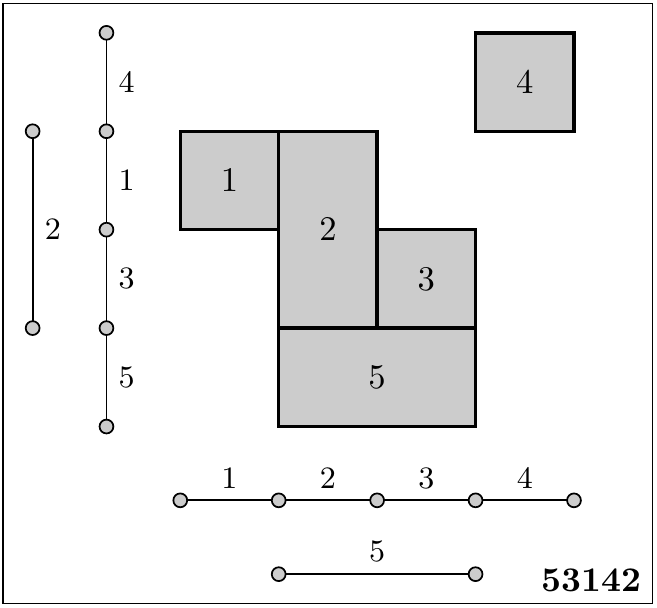}\\ 
 \includegraphics[width=\wsp\textwidth]{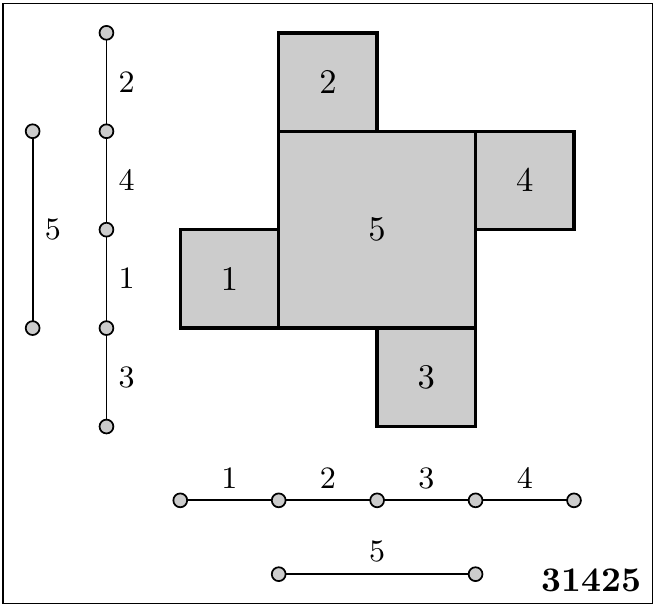} &
\multicolumn{1}{c}{} & 
\multicolumn{1}{c}{}\\ 
\end{tabular}
\caption{Cliques solving the construction problem for given $\Gamma_1$ and $\Gamma_2.$\protect\label{piatki}}
\end{figure}

 Now, observe that each clique $\ka C$ is determined by an isomorphism  of $G_2$ onto $\Gamma_2$ and vice versa.   These isomorphisms can be expressed in terms of appropriate labellings of $\Gamma_2$: We write  the cycle $G_2$ as a  sequence $l_1l_2l_3l_4l_5$ and define the corresponding labelling
$$  
\text{$[a_{2,k},a_{2,k+1}]\mapsto l_{k+1}$ for $k=0,\ldots, 3$ and $[a_{1,1},a_{1,3}]\mapsto l_5$.} 
$$
We have ten sequences corresponding to $G_2$: $13524$, $35241$, $52413$, $24135$, $41352$, $14253$, $42531$, $25314$, $53142$, $31425$.  Thus, we have ten cliques, which are depicted in Figure \ref{piatki}. Each clique is labelled by an appropriate  sequence.
\begin{uwa}
From a geometric point of view,  we have  essentially two different types of 5-cliques in $\ka I^2$. They are exemplified  by cliques with labels $24135$ and $35241$.  Both types were previously  mentioned in \cite[Figure 27]{Simon1}.
\end{uwa}

\section{Maximum incompressible cliques in \texorpdfstring{$\ka I^3$}{dimension 3} \label{mincl}}

Our present goal is to describe all cliques of maximum cardinality in $\ka I^3$. We begin with a general discussion concerning cliques in $\ka I^n$. 

Let us suppose that $\ka G\subset \ka I^n$ is a clique of cardinality  $b_n$, where $n\ge 2$. Let $\ka G_i=\{ J_i\colon J \in \ka G\}$, $i\in [n]$.  Let $s_i=s(\ka G_i)$ (see: Proposition \ref{shom})  and $f_i$ be a homomorphism of $\ka G_i$ into $\ka I(s_i)$.  Let us set  $f=f_1\times\cdots \times f_n$  and
$$
\ka C= f(\ka G)=\{f_1(J_1)\times \cdots \times f_n(J_n)\colon J\in \ka G\}.
$$ 
Clearly, $\ka C$ is a clique of maximum cardinality $b_n$ in $\ka I^n$. 
Let $\gamma_i =\sum_{I\in \ka C} I_i$.  By the definition of $\ka C$,   the formal sum $\gamma_i\in \enn\ka I(s_i)$ is incompressible.  Therefore,  by  Propositions \ref{irred},  and  \ref{oszac},
\begin{equation}
\label{esi}
s_i\le 2\alpha(\gamma_i)-1\le 2b_{n-1}-1.
\end{equation}
And if $s_i\ge 3$, then again by Proposition \ref{irred},
\begin{equation}
\label{nosnik}
\ka S(s_i)\subseteq \supp \gamma_i.
\end{equation}

Let us remark that  a description  of cliques like $\ka C$, call them \textit{incompressible}, is a strictly combinatorial problem: All such cliques consist of intervals $I$ contained in  $[0, 2b_{n-1}-1]^n$ whose vertices have integer coordinates.  Suppose that this problem can be effectively solved. As we have observed, if $\ka G$ is a clique of maximum cardinality, then  there is a mapping $f$ which sends it to an incompressible clique. We can hope that this fact can be employed in order to give a  full description of all maximum cliques in $\ka I^n$.   As we shall see, this general strategy works for $n=3$.

\begin{npr}
\label{np12klika}
Let $\ka C\subset \ka I^3$ consists of the following intervals:
$$
\begin{array}{rc}
1.\,\, & [0,1]\times[3,4]\times[1,4] \\
2.\,\, & [0,1]\times[4,5]\times[2,4] \\  
3.\,\, & [1,2]\times[3,4]\times[1,3] \\
4.\,\, & [1,2]\times[4,5]\times[2,3] \\
5.\,\, & [1,3]\times[1,2]\times[3,4] \\
6.\,\, & [1,4]\times[0,1]\times[3,4] \\ 
7.\,\, & [2,3]\times[1,2]\times[4,5] \\
8.\,\, & [2,4]\times[0,1]\times[4,5] \\
9.\,\, & [3,4]\times[1,3]\times[1,2] \\
10.\,\,&[3,4]\times[1,4]\times[0,1] \\
11.\,\,&[4,5]\times[2,3]\times[1,2] \\
12.\,\,&[4,5]\times[2,4]\times[0,1]  
\end{array}.
$$
One can easily check that $\ka C$ is an incompressible clique in $\ka I^3$; that is, each $\gamma_i=\sum_{I\in \ka C} I_i$ is incompressible  (in fact, $\gamma_1=\gamma_2=\gamma_3$). As $|{\ka C}|=12$, it follows that $b_3\ge 12$.  

\begin{figure}[ht]
	\centering
		\includegraphics[width=0.7\textwidth]{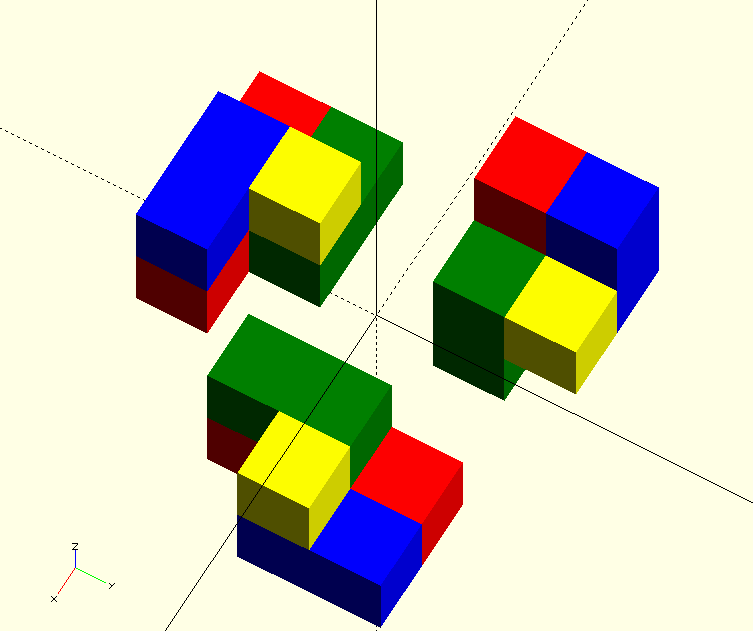}
	\caption{An illustration of the clique $\ka C$ described in Example 1. The clique is translated so that the circumscribed box  is centred at the origin.}
	\label{fig:example1}
\end{figure}

\end{npr}

Now, we distinguish  a set $L$  of combinations $\lambda\in \enn\ka I$ such that for every incompressible clique $\ka C$ in $\ka I^3$ of the  maximal cardinality,  each $\gamma_i=\sum_{I\in \ka C} I_i$  is a  member of $L$.  This set is characterized by the following conditions:
 \begin{description}
\item[(A)] $12\le| \lambda|\le 17$;
\item[(B)] $\alpha(\lambda)\le 5$;
\item[(C)]  there is $s$, $3\le s\le 9$, such that $\ka S(s)\subseteq \supp \lambda\subseteq \ka I(s)$.
\end{description}

Condition (A) reflects the fact that $ 12\le b_3\le 17$ (the second inequality is derived from Proposition \ref{rekurencja}) ; (B) is implied by  Proposition \ref{oszac} and the equality $b_2=5$.  Condition (C) reflects (\ref{esi}) and (\ref{nosnik}). Only the assumption $s\ge 3$ requires clarification. Suppose $\lambda= \gamma_1$, where    
$\gamma_1=\sum_{I\in \ka C} I_1$ and $\ka C$ is an incompressible clique of cardinality $b_3$ in $\ka I^3$.  If  there would be $s=s_1\le 2$, then  
$$
12\le |\ka C|=|\gamma_1|\le 2\alpha(\gamma_1)\le 10,
$$
which is impossible.

If one wants to determine all members of $L$, then  one needs an effective method  to verify (B).  We discuss this question now. 

Let $s$  be a positive integer. For $\varepsilon \in \{0,1\}^s$, let the subfamily $\ka I(s,\varepsilon)$ of $\ka I(s)$ be defined as follows
\begin{quote}
$[u,v]\in \ka I(s,\varepsilon)$ if and only if   $1\le u < v\le s$ and $\varepsilon_u=0, \varepsilon_v=1$ or $v=1$ and $\varepsilon_1=1$ or $u=s$ and $\varepsilon_s=0.$  
\end{quote}
It is clear that all the families $\ka I(s, \varepsilon)$, $\varepsilon\in \{ 0,1\}^s$, are independent sets with respect to the adjacency defined in $\ka I^1$. Moreover, each maximal independent subset of $\ka I(s)$ can be found among these families.  This observation leads to the following
\begin{pr}
Let $s$ be a positive integer.  If $\lambda=\sum k_I I$ belongs $\enn\ka I(s)$, then   
$$
\alpha(\lambda)=\max \{ \sum_{I\in \ka I(s,\varepsilon)} k_I\colon \varepsilon \in \{0,1\}^s\}
$$
\end{pr}
This proposition enables us to identify the problem of determining $L$ as a problem of integer linear programming.

 Let 
$ L(s,v)$ be the set of all these  $\lambda\in L$ for which $\ka S(s)\subseteq \supp \lambda\subseteq \ka I(s)$ (comp. (C)) and $|\lambda|=v$. Clearly,   
$L$ is a union of the sets $L(s,v)$, $3\le s\le 9$, $12\le v \le 17$.  A combination $\lambda=\sum_{I\in \ka I(s)} k_II$ is a member of $L(s,v)$ if and only if $k_I$, $I\in \ka I(s)$, satisfy the following system of linear inequalities:
\begin{eqnarray}
\sum_{I\in \ka I(s)}k_I &= & v;\\
\sum_{I\in \ka I(s,\varepsilon)} k_I &\le &5, \, \text{for every $\varepsilon\in \{0,1\}^s$;}\\
k_I & \ge & 1,\, \text{for $I\in \ka S(s)$.}
\end{eqnarray}

All solutions of  this system can be found using standard mathematical packages. We have preferred to run a simple Python code, and to perform tests using  SCIP, a mixed integer programming solver.  As we have computed, $L(s,14)$ is non-empty only for $s=9$.  Moreover, $L(9,14)$ consists of a single element   
$$
\bar{\lambda}= [0,1]+ [1,2]+[2,3]+[3,4]+[4,5]+[5,6]+[6,7]+[7,8]+[8,9]+[9,10]+[1,3]+[3,5]+[5,7] +[7,9].  
$$
The sets $L(s,v)$, for $v>14$, appear to be empty, which can also be easily deduced from the uniqueness of $\bar{\lambda}$.   

One can  check by hand  that  $E(\bar{\lambda})=20$.  The following proposition is  an immediate consequence of the definition of adjacency in $\ka I^n$. 
\begin{pr}
\label{stopnie}
Let $\ka C$ be a subset of  $\ka I^n$ and  $\gamma_{i}=\sum_{I\in \ka C}  I_i$, $i \in [n]$, then
$$  
|E(\ka C)|\le \sum_{i=1}^n |E(\gamma_i)|.
$$ 
\end{pr}
Suppose that there is a clique in $\ka I^3$ which has 14 elements. Then there would be an incompressible clique $\ka C$ of the same cardinality.  By the uniqueness of $\bar{\lambda}$, we would have   $\bar{\lambda}=\sum_{I \in \ka C} I_i$ for $i\in [3]$, and  by Proposition \ref{stopnie}
$$
91=|E(\ka C)|\le 3|E(\bar{\lambda})|=60.
$$ 
This contradiction  implies 
$$
b_3\le 13.
$$
It also shows that from the point of view of maximum cliques only the cases $v=12, 13$ are of interest. The cardinalities of the sets $L(s,v)$, $3\le s \le 9$, $12\le  v \le 13$,  are collected in Table \ref{licznosc}. 

 \begin{table}
\label{licznosc}
\begin{center}
\begin{tabular}{|l|*{7}{c|}}\hline
\diagbox{$v$}{$s$} & 3&4&5&6&7&8&9\\
\hline
12 &5&19&210&164&82&13& 1\\
13 &0&  0&  20&  35&55&13& 3\\\hline
\end{tabular}
\end{center}
\caption{Cardinality of $L(s,v)$. }
\end{table} 

 Let $\aut s$ be the automorphism group of  the graph $\ka I(s)$. By the correspondence $\aut s \ni f  \mapsto f_*$, it acts on $\enn \ka I(s)$ in a natural way. Therefore, for every $L(s,v)$, one can form the quotient set $L(s,v)/\aut s$. Let $\hat{L}(s,v)$ be a system of representatives (a selection) of   $L(s,v)/\aut s $.  Let $L(v)$ be the union of all $L(s,v)$, $3\le s\le 9$, for $v=12,13$, and $\hat L(v)$ be the union of the corresponding sets $\hat L(s,v)$. 
It should be clear that any maximum incompressible clique $\ka D$ of $\ka I^3$ can be recovered (by an automorphism) from a clique $\ka C$ such that $\gamma_i=\sum_{I \in \ka C} I_i\in \hat L(v)$, for $i\in [3]$, where $v=|\ka D|$. Therefore, we may always consider  $\hat L(v)$ instead  of $L(v)$.  This will slightly simplify further considerations, as in general the sets $\hat L(s,v)$ are smaller than $L(s,v)$.  
In order to fix $\hat L(v)$ we work with, we collect all the sets $\hat L(s,v)$ in Appendix B. We present here  only the table containing the cardinalities of $\hat L(s,v)$ to allow  the reader to compare it with Table~\ref{licznosc}. 
\begin{table}
\label{licznosc2}
\begin{center}
\begin{tabular}{|l|*{7}{c|}}\hline
\diagbox{$v$}{$s$} & 3&4&5&6&7&8&9\\
\hline
12 &1&4&37&29&21&5& 1\\
13 &0&  0&  5&  8&13&5& 2\\\hline
\end{tabular}
\end{center}
\caption{Cardinality of $\hat L(s,v)$. }
\end{table}      
 
We understand now that it makes sense to solve the following `restricted' construction problem as a step towards the classification of all maximum cliques in $\ka I^3$: 
\begin{quote}
For every triple $  \vv{\gamma}=(\gamma_1,\gamma_2,\gamma_3)$ in $\hat L(v)^3$, where $v=12,13$, find all cliques $\ka C$ in $\ka I^3$ such that
$\gamma_i=\sum_{I\in \ka C} I_i$  for every $i\in [3]$.   
\end{quote}

Let us define the sets $N(v)$ for $v=12,13$:
$$
N(v)=\left\{\vv{\gamma} \in \hat L(v)^3\colon E(\gamma_1)+E(\gamma_2)+E(\gamma_3)\ge \binom{v}{2}\right\}.
$$
Proposition \ref{stopnie} shows that if the construction problem  has a solution for  $\vv{\gamma}\in\hat L(v)^3$, then $\vv{\gamma}\in N(v)$.  Therefore,  we can consider our problem only for elements of the  sets $N(v)$. 
We can make further restrictions by selecting a triple from each class of triples equivalent up to reordering of components.  However, we want to make our choice somewhat special. 

We may distinguish between two types of $\lambda\in \hat L(v)$:
\begin{description}
\item{(I)}  There are induced combinations  $\beta^1$ and $\beta^2$ of $\lambda$ (see Introduction for the definition) such that
\begin{description}  
\item{(I.\,1)} $\quad|\beta^j|=\alpha(\beta^j)=5$,  for every $j\in [2]$, 
\item{(I.\,2)} $\quad\supp \beta^1\cap\supp \beta^2=\emptyset$.
\end{description}
\item{(II)}  $\lambda$ is not of type I.
\end{description}
Let us remark that  (I.\,1) means that  $G_{\beta^j}$ are 5-anticliques.  As it can be easily computed (see: Appendix B), 
each $\lambda$ of type II satisfies the inequality 
$$
E(\lambda)\le \frac 1 3 \binom{v}{ 2} .
$$
Moreover, this inequality is sharp with only one exception; that is,
$$
\lambda^*=[0, 1]+[1, 2]+[1, 3]+ [1, 5]+[2, 3]+[2, 5]+[3, 4]+2[3, 5]+[4, 5]+2[5, 6].
$$
(Clearly, $\lambda^*$ belongs to $\hat L(5,12)$ then).  As a consequence, 
there is only one $\vv{\gamma}\in N(v)$ such that 
all $\gamma_i$ are of type II, in which case they are equal to $\lambda^*$.   The following result shows that  this particular $\vv{\gamma}$ can be excluded from further considerations. 
\begin{pr}
\label{bezrozw}  
The construction problem has no solutions for $\vv{\gamma}=(\lambda^*, \lambda^*, \lambda^*)$.
\end{pr}
\proof 
Conversely, suppose there is a clique $\ka C$ in $\ka I^3$  such that $\lambda^*=\sum_{I\in \ka C} I_i$, whenever $i\in [3]$.  As $[5,6]\in \supp \lambda^*$, there is $K \in \ka C$  for which $K_1=[5,6]$.  Similarly, there are $I^k\in \ka C$  such that $I^k_1=[k,5]$, for $k\in [4]$. Since $[3,5]$ occurs with multiplicity $2$ in $\lambda^*$, there is $J^3 \in \ka C$  which is different from $I^3$, and satisfies the equation $J^3_1=[3,5]$.  Since $|E(\lambda^*)|=22$ and $|E(\ka C)|=66$, it follows from the definition of adjacency in $\ka I^3$ and Proposition \ref {stopnie}  that for every pair $C,D \in \ka C$  there is only one $i\in [3]$ for which the intervals $C_i$ and $D_i$ are adjacent (i.e. intersect at a single point).  Since all $I^k_1$, $k\in [4]$, form an anticlique, we deduce that $\ka B=\{I^k_{\{1\}^\mathsf{c}}\colon k\in[4]\} \cup \{J^3_{\{1\}^\mathsf{c}}\}$ is a 5-clique in $\ka I^2$.  Therefore, as is shown in Section \ref{maxc 2D}, the families $\ka B_i=\{I^k_{i}\colon k\in[4]\} \cup \{J^3_{i}\}$, $i=2,3$, have to be 5-cycles contained in $\supp \lambda^*$.  Let us collect all  possible 5-cycles with vertices in $\supp \lambda^*$: 
\begin{eqnarray}
1. &&  [0, 1], [1, 2], [2, 3], [3, 4], [1, 3]\nonumber \\
2. &&  [0, 1], [1, 2], [2, 5], [5, 6], [1, 5]\nonumber\\
3. &&  [2, 3], [3, 4], [4, 5], [5, 6], [3, 5]\nonumber\\
4. &&  [0, 1], [1, 3], [3, 5], [5, 6], [1, 5]\nonumber\\
5. &&  [1, 2], [2, 3], [3, 5], [5, 6], [2, 5]\nonumber
\end{eqnarray}
As $K_1$ is adjacent to each of the intervals $I^k_1$, $k\in[4]$,  interval $K_i$ cannot be adjacent to any of the intervals belonging to $\ka B_i$, for $i=2,3$.  Consequently,  if  a 5-cycle among listed were equal to $\ka B_i$, then there would be an interval  in $\supp\lambda^*$ which is not adjacent to any member of  this cycle. Cycles no.~4 and no.~5 do not conform this condition,  and as such can be eliminated. On the other hand, for cycle no.~1 the only existing interval is $[5,6]$, for  cycle no.~2, the interval $[3,4]$ and for cycle no.~3, the interval $[0,1]$. Thus, $\{K_2, K_3\}\subset \{[0,1], [3,4], [5,6]\}$. Let us remind that $[5,6]$, similarly as $[3,5]$, occurs in $\lambda^*$ with multiplicity 2. Therefore, there is $L\in \ka C\setminus\{K\}$ such that $L_1=[5,6]$. By the same argument as applied to $K$, $\{L_2, L_3\}\subset \{[0,1], [3,4], [5,6]\}$.  Since  $[0,1], [3,4], [5,6]$ form an anticlique, $K$ and $L$ are not adjacent, which contradicts the assumption  that $\ka C$ is a clique. \hfill$\square$                

Therefore, by Proposition \ref{bezrozw}, we may assume that at least one of the components of $\vv{\gamma}\in N(v)$ is of type I . We may declare that $\gamma_3$ is such a component. Suppose there is a clique $\ka C \subset \ka I^3$ of cardinality $v$ to be guessed such that $\gamma_i=\sum_{I\in \ka C} I_i$.  Let $\beta^j$, $j\in [2]$, be combinations induced from $\gamma_3$, as  described in the definition of type I.  Since these combinations have disjoint  supports there are disjoint subfamilies $\ka C^j$, $j\in [2]$, of $\ka C$  such that $\beta^j=\sum_{I\in \ka C^j} I_3$.   Since $\beta^j$ are 5-anticliques,  $\ka C^j_{\{1,2\}}=\{I_{\{1,2\}}\colon I \in \ka C^j\}$  have to be   5-cliques in $\ka I^2$.  Thus,  the families $\ka C^j_i=\{I_i\colon I \in \ka C^j\}$,  $i,  j \in [2]$,  are 5-cycles, as is explained in Section \ref{maxc 2D}.  

 Let us write
$$
\gamma=\sum_{I \in \ka C} I,\quad \gamma^j= \sum_{I \in \ka C^j} I.
$$
Clearly, $\gamma\ge \gamma^1+\gamma^2$.  Therefore, there is $\gamma^3\in \enn\ka I^3$ such that
$
\gamma=\sum_{j=1}^3\gamma^j\, \text{and}\,  |\gamma^3|=v-10. 
$
Consequently,  
\begin{equation} 
\label{rozklad}
\gamma_{\{1,2\}}=\sum_{j=1}^3 \gamma^j_{\{1,2\}}\quad  \text{and}\quad \gamma_i=\sum_{j=1}^3\gamma^j_i\quad  \text{for $i\in [2]$,}  
\end{equation}
where each $\gamma^j_i$, $j\in[2]$, is the formal sum of all intervals constituting the 5-cycle $\ka C^j_i$. Moreover, by Corollary  \ref {wnor}   
\begin{equation}
\label{prostokaty}
\alpha(\gamma_{\{1,2\}}) \le \omega(\gamma_3)\le \omega(\ka I^1)\le 2.
\end{equation}

Now, we are prepared to establish a \textbf{procedure for finding all incompressible cliques of maximum cardinality in} $\ka I^3$.  Since each $\lambda\in \hat  L(v)$ can potentially be  equal to $\gamma_1$ or $\gamma_2$ for a certain maximum clique in $\ka C$, we have to produce  all possible decompositions $\lambda=\sum_{j=1}^3 \lambda^j$, as described in the second part of (\ref{rozklad}).  In other words, we have to extract all possible quadruples $q=(\lambda^1,\lambda^2, \lambda^3, \lambda)$ such that $\lambda^j$ are 5-cycles for $j\in [2]$ and $\sum_{j=1}^3\lambda^j$ equals $\lambda$.  To this end, we need to determine all  5-cycles with their supports contained in  $\ka I(9)$. 
\skok
\textbf{Step 1.}  Find the set $\text{Co}_5$  
consisting of all 5-cycles contained in $\enn\ka I(9)$.  
\skok
Let us recall (see: Section \ref{maxc 2D}) that  there is a one-to-one correspondence between the sets of all 5-cycles in  $\enn\ka I(9)$, and in $\ka I(9)$ given by the mapping $\kappa\mapsto \supp \kappa$.   From a technical point of view, the latter set is calculated rather than $\text{Co}_5$.  

\begin{figure}
\centering 
\includegraphics{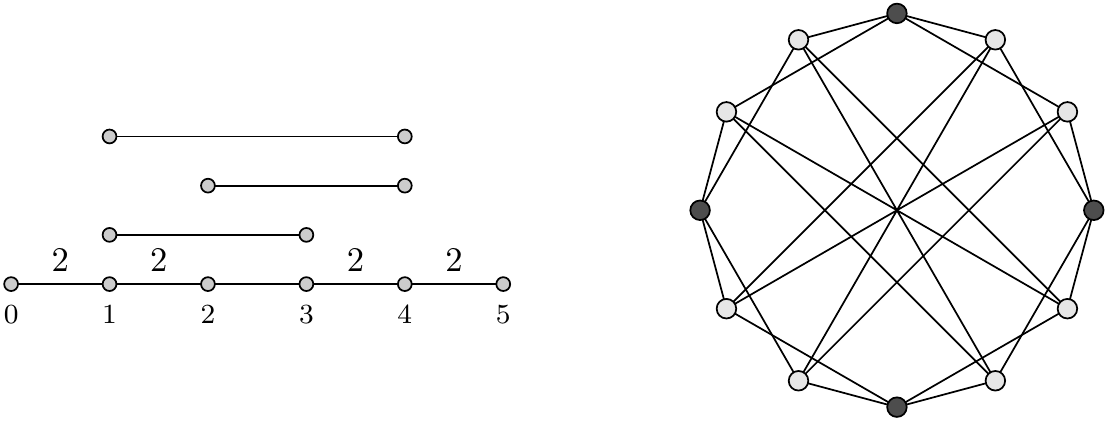}
\caption{\label{trefle} Diagrams for $\lambda^{\clubsuit} $.}
\end{figure}

\begin{figure}
\centering
\includegraphics{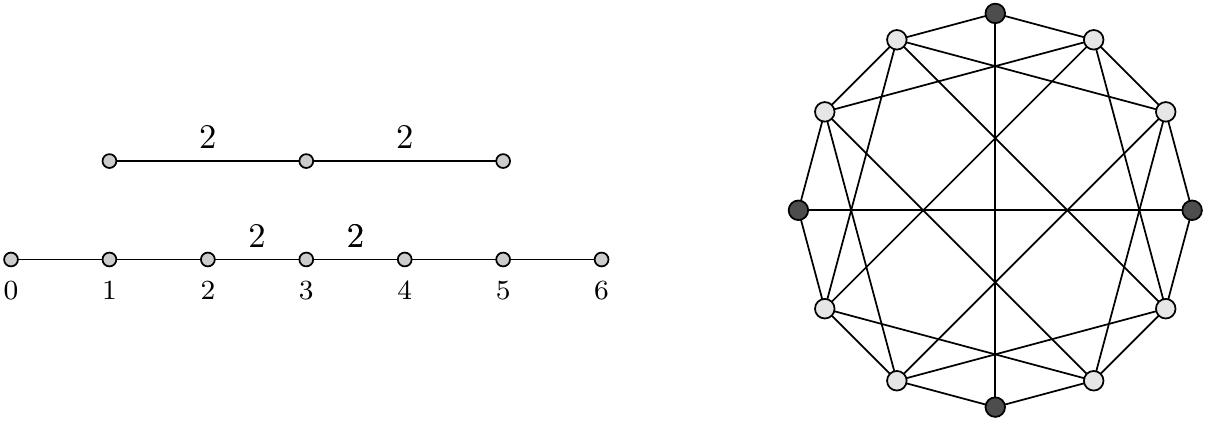}
\caption{Diagrams for $\lambda^{\spadesuit}$. \protect\label{piki} }
\end{figure}

\begin{figure}[ht]
\centering
\includegraphics{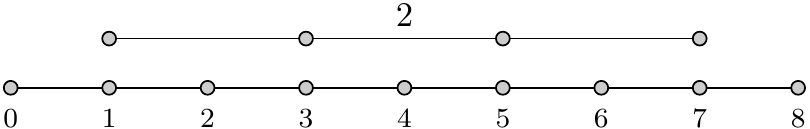}
\caption{Diagram for $\lambda^{\diamondsuit}$.\protect \label{cara}} 
\end{figure}

The first two components of the quadruples  under consideration are 5-cycles, therefore, it seems reasonable to determine the Cartesian product $\text{Co}_5\times \text{Co}_5$. Since $| \text{Co}_5|=118$, the latter set has 13924 elements. Observe however that not all pairs of 5-cycles can be components of a quadruple. For example,  if $\Gamma=[0,1]+[1,2]+[2,3]+[3,4]+[1,3]$ and  $\Delta=[0,1]+[1,3]+[3,7]+[7,8]+[1,7]$ would be such components for $v=12$, then there would exist $\lambda\in \hat L(12)$ such that $\Gamma+\Delta\le \lambda$. If we take into account  (C),  the support of such a $\lambda$ has to contain one of the two sets  $A=\ka S(8)$ or $A=\ka S(9)$, which is impossible as $A-(\supp(\Gamma)\cup\supp(\Delta))$ has more than two elements. To define the set of appropriate pairs of 5-cycles, we need the function $\Gamma\mapsto\operatorname{sep} (\Gamma)$ which determines the smallest $s$ such that $\supp(\Gamma)\subset \ka I(s)$. As $\Gamma$ can be written as a combination $[a_0,a_1] + [a_1,a_2] + [a_2,a_3] + [a_3,a_4] + [a_1,a_3]$,  where $a_0<\ldots<a_4$, and all $[a_i,a_j]$ are in $\ka I(9)$, we have
$$
 \operatorname{sep}\, (\Gamma)=
\left\{
\begin{array} {ll}
 a_3, &\text{if $a_4=a_3+1$,}\\
a_4,&\text{otherwise}.\\
\end{array}
\right.
$$ 
\skok 
\textbf{Step 2.}   For $v=12,13$, find the set $\text{Cp}_5(v)$, consisting of all pairs $(\Gamma, \Delta)\in \text{Co}_5\times \text{Co}_5 $ such that $|\ka S(s)-(\supp(\Gamma)\cup\supp(\Delta))|\le v-10$, where $s=\max\{\operatorname{sep}\,(\Gamma), \operatorname{sep}\,(\Delta)\}$. 
\skok
The sets $\text{Cp}_5(v)$ are much smaller compared with $\text{Co}_5\times \text{Co}_5$, and yet $(\lambda^1, \lambda^2)$ belongs to $\text{Cp}_5(v)$ for every admissible quadruple $q=(\lambda^1, \lambda^2, \lambda^3, \lambda)$ .  We have   $|\text{Cp}_5(12 )|=384 $  and $|\text{Cp}_5(13)|=1135$.  
\skok
\textbf{Step 3.} For $v=12, 13$, find the set 
$\text{Mc}(v)$ of all quadruples $(\lambda^1, \lambda^2, \lambda^3, \lambda)$ such that  $(\lambda^1, \lambda^2)\in \text{Cp}_5(v)$, $\lambda\in \hat L(v)$,   $\lambda^1+ \lambda^2\le \lambda$ and  $\lambda^3=\lambda-\lambda^1-\lambda^2$. 
\skok 
The cardinalities of the sets $\text{Mc}(v)$ are as follows:  $|\text{Mc}(12)|= 372$, $|\text{Mc}(13)|=409$.   The next step is to construct all candidates   for $\gamma_{\{1,2\}}$. We want  them to be decomposed as described in (\ref{rozklad}) and  satisfy (\ref{prostokaty}). To this end, for every pair $q= (\lambda^1, \lambda^2, \lambda^3, \lambda)$, $\bar{q}=(\bar\lambda^1, \bar \lambda^2, \bar\lambda^3, \bar\lambda)$ belonging to $\text{Mc}(v)$, let us define the set $\text{Fl}(q, \bar q)$  of all quadruples $p=(\varphi^1,\varphi^2, \varphi^3, \varphi)\in (\enn\ka I^2)^4$ such that 
\begin{enumerate}
\item[(F1)]
$G_{\varphi^j}$ are 5-cliques in $\ka I^2$, for $j\in [2]$, 
\item[(F2)]
$\varphi = \varphi^1 + \varphi^2 + \varphi^3$, 
\item[(F3)]
if $p_i= (\varphi^1_i,\varphi^2_i, \varphi^3_i, \varphi_i)$, $i\in[2]$, then $p_1=q$, $p_2=\bar q$;
\item[(F4)]
$\alpha(\varphi)\le 2$.
\end{enumerate}
As is stated in Sect. \ref{maxc 2D}, there are 10 ways to combine two $5$-cycles   from $\ka I$  to get a $5$-clique in $\ka I^2$. Moreover, there are at most $(v-10)!$ ways to merge $\lambda^3$ with $\bar \lambda^3$ in order to get $\varphi^3$. Therefore,  the cardinality of $\text{Fl}(q,\bar q)$ is at most $100(v-10)!\le 600$.  Condition  (F4) can cut the latter number substantially. 
\skok
\textbf{Step 4.} Find the sets $\text{Flat}(v)$, $v=12, 13$, which  are the unions of  $\text{Fl}(q,\bar q)$, when $(q,\bar q)$ runs over $\text{Mc}(v)^2$.  
\skok
It appears  that $\text{Flat}(13)$ is empty. Consequently, \textbf{no 13-cliques in $\ka I^3$ exist.} By Example \ref{np12klika}, we conclude:
\begin{tw}
\label{be3}
$$
b_3=12.
$$
\end{tw} 
$\text{Flat}(12)$ contains 96 elements.  Interestingly enough,  if $(\varphi^1,\varphi^2, \varphi^3, \varphi)$ belongs  to  $\text{Flat}(12)$, then $\varphi_1=\varphi_2$. Moreover, $\varphi_1$ is equal to one of the  three combinations:
\begin{eqnarray}
\lambda^\clubsuit & = & 2[0,1]+2[1,2]+[1,3]+[1,4]+[2,3]+[2,4]+2[3,4]+2[4,5],\nonumber \\
\lambda^\spadesuit & = & [0,1]+[1,2]+2[1,3]+2[2,3]+2[3,4]+2[3,5]+[4,5]+[5,6],\nonumber \\
\lambda^\diamondsuit & = & [0,1]+[1,2]+[1,3]+[2,3]+[3,4]+2[3,5] +[4,5]+[5,6]+[5,7]+[6,7]+[7,8].\nonumber
\end{eqnarray}  
Figures \ref{trefle}--\ref{cara} show self-explanatory  diagrams for them.  Combinations $\lambda^\clubsuit$ and  $\lambda^\spadesuit$ are of type I, as exemplified by Table \ref{decomp}, while   $\lambda^\diamondsuit$ is of type II.

\begin{table}[hbp]
\begin{center}
\begin{tabular}{|c|c|c|}
\hline
& $\beta^1$ & $\beta^2$\\
\hline
$\lambda^\clubsuit$ & 2[0,1]+[2,4]+2[3,4] & 2[1,2]+[1,3]+2[4,5] \\

$\lambda^\spadesuit$ & [0,1]+2[3,4]+2[3,5] & 2[1,3]+2[2,3]+[5,6]\\
\hline
\end{tabular}
\end{center}
\caption{$\lambda^\clubsuit$ and $\lambda^\spadesuit$ are of type I. \protect\label{decomp}}
\end{table}
\begin{pr}
\label{marozw}
If the construction problem has a solution for $\vv{\gamma}=(\gamma_1,\gamma_2,\gamma_3) \in N(12)$, then  
$$\vv{\gamma}=(\lambda^\clubsuit, \lambda^\clubsuit, \lambda^\clubsuit) \quad \text{or} \quad \vv{\gamma}=(\lambda^\spadesuit, \lambda^\spadesuit, \lambda^\spadesuit).
$$
\end{pr}

\proof  Let us remind that we may assume $\gamma_3$ to be of type I.  Then it follows from the definition, and properties of $\text{Flat}(12)$  that  $\gamma_1=\gamma_2$ and $\gamma_1\in \{\lambda^\clubsuit,\lambda^\spadesuit, \lambda^\diamondsuit\}$.  

Suppose first that $\gamma_1=\lambda^\diamondsuit$.  Recall that $|E(\lambda^\diamondsuit)|=19$.  Since $\vv{\gamma}\in N(12)$, it can be deduced that $|E(\gamma_3)|\ge 28$.  There is a unique element in $\hat L(12)$ which satisfies this inequality. Namely, $\gamma_3= 2[0,1] +2[1,2] +3[1,3]+2[2,3]+3[3,4]$.  ($|E(\gamma_3)|=29$, to be precise).  Suppose that there is a 12-clique $\ka C$ such that  $\gamma_i=\sum_{I\in \ka C}I_i$, $i\in [3]$.   Let $\beta^1=2[0,1]+3[3,4]$, $\beta^2=2[1,2]+3[3,4]$ and $\beta^3=3[1,3]+2[2,3]$.  As  all $\beta^j$ are subcombinations of  $\gamma_3 $,  there are  $\gamma^j\le \gamma=\sum_{I \in \ka C} I$ for which $\gamma_3^j=\beta^j$.  Observe now that $\alpha(\beta^j)=5$ for each $\beta^j$. We already know that it implies  each $\gamma^j_1$  to be a 5-cycle.  Since  combinations $\beta^1$, $\beta^3$ have disjoint supports, we deduce that $\gamma^1_1+\gamma^3_1\le \gamma_1=\lambda^\diamondsuit$. There are only two 5-cycles 
such that their sum does not exceed $\lambda^\diamondsuit$; these are $\kappa=[0,1]+[1,2]+[2,3]+[3,4]+[1,3]$  and $\kappa'=[4,5]+[5,6]+[6,7]+[7,8]+[5,7]$. Thus, $\{\gamma_1^1, \gamma_1^3\}=\{\kappa,\kappa'\}$. Similarly,  by the fact that    $\beta^2$, $\beta^3$ have disjoint supports, we obtain $\{\gamma_1^2, \gamma_1^3\}=\{\kappa,\kappa'\}$.  Since $\gamma_1=\lambda^\diamondsuit$, there is $J\in \ka C$  such that $J_1=[3,5]$.  Since the supports of $\beta^j$, $j\in [3]$ form a  covering of $\supp\gamma_3$, there is $k$  such that $J_3\in \supp\beta^k$.  As $\beta^k$ is an induced subcombination of $\gamma_3$, we deduce that $J\in \supp \gamma^k$. Therefore,
$$
[3,5]\in \supp \gamma^k_1\subset\supp \kappa \cup \supp\kappa' =\{[i,i+1]\colon i =0,\ldots, 7\} \cup\{[1,3], [5,7]\}, 
$$
which is a contradiction.

If $\gamma_1\in \{\lambda^\clubsuit, \lambda^\spadesuit \} $, then it is of type I. By symmetry, it  enforces $\gamma_2=\gamma_3$. Thus $\vv \gamma\in \{(\lambda^\clubsuit, \lambda^\clubsuit, \lambda^\clubsuit)$, $(\lambda^\spadesuit, \lambda^\spadesuit, \lambda^\spadesuit)\}$. \hfill $\square$

Now, the final steps of  our construction  follow.    Let us extract  two sets from \text{Flat}(12):
$$
\begin{array}{l}
\text{Fl}^\clubsuit  =  \{ (\varphi^1,\varphi^2, \varphi^3)\colon  (\varphi^1,\varphi^2, \varphi^3, \varphi)\in \text{Flat}(12), \,\, \varphi_1=\lambda^\clubsuit\} ,\\ 
\\
\text{Fl}^\spadesuit  = \{ (\varphi^1,\varphi^2, \varphi^3)\colon  (\varphi^1,\varphi^2, \varphi^3, \varphi)\in \text{Flat}(12), \,\, \varphi_1=\lambda^\spadesuit\}.
\end{array}
$$
Let  $(\beta^1, \beta^2, \beta^3)$ be defined so that $\lambda^\clubsuit=\beta^1+\beta^2 +\beta^3$ where  $\beta^1$, $\beta^2$ are the elements  of $\lambda^\clubsuit$-row of Table \ref{decomp} .   
 \skok
\textbf{Step 5} Find 
$$
\text{Cq}^\clubsuit= \{\pi^1+\pi^2 +\pi^3\in \enn\ka I^3\colon   (\pi^1_1,\pi^2_1, \pi^3_1)\in \text{Fl}^\clubsuit, \,   (\pi^1_2,\pi^2_2, \pi^3_2)=(\beta^1, \beta^2, \beta^3), \,\, \alpha(\pi^1+\pi^2 +\pi^3)=1\}.
$$
\skok
Therefore, $\text{Cq}^\clubsuit$ consists of all  12-cliques  $\ka C\subset \ka I^3$, written as formal combinations, such that $\sum_{I\in \ka C} I_i=\lambda^\clubsuit$ for $i \in [3]$.  It appears that we have 
\begin{equation}
\label{ile trefli}
|\text{Cq}^\clubsuit|=64.
\end{equation}

In the same manner, we define  $\text{Cq}^\spadesuit$.
\skok
\textbf{Step 6} Find   $\text{Cq}^\spadesuit$.
\skok
The elements of $\text{Cq}^\spadesuit$  correspond to 12-cliques  $\ka C\subset \ka I^3$ such that    $\sum_{I\in \ka C} I_i=\lambda^\spadesuit$. We have  
\begin{equation}
\label{ile piki}
|\text{Cq}^\spadesuit|= 256.
\end{equation} 
\begin{tw}
\label{allcin}  $\ka C$ is an incompressible maximum clique in $\ka I^3$ if and only if  $\,\sum_{I\in \ka C}I \in \operatorname{Cq}^\clubsuit\cup\operatorname{Cq}^\spadesuit$.
\end{tw}
\proof
$(\Rightarrow)$
Let  $\gamma=\sum_{I\in \ka C} I$. By  Theorem \ref{be3} and  the definition of the sets $L(s,v)$, there are $s_i$, $i \in [3]$, such  that $\gamma_i\in L(s_i,12)$. Let   automorphisms  $f_i\in  \aut {s_i}$  be  chosen so that  $({f_i})_*(\gamma_i)\in \hat L(s_i,12)$.  Then  $( ({f_1})_*(\gamma_1),({f_2})_*(\gamma_2), ({f_3})_*(\gamma_3))$ has to belong   to $N(12)$.   By   Proposition  \ref{marozw},  $({f_i})_*(\gamma_i)=\lambda^\clubsuit$ for $i\in [3]$ or $({f_i})_*(\gamma_i)=\lambda^\spadesuit$ for $i\in [3]$. Observe now that  $\lambda^\clubsuit$  is  invariant under the action of $\aut 4 $; that is, if  $g\in \aut 4$, then $g_*( \lambda^\clubsuit)=\lambda^\clubsuit$. Similarly,   $\lambda^\spadesuit$ is  invariant under the  action of $\aut 5$. Therefore, $(\gamma_1, \gamma_2, \gamma_3)$ is one of the triples  $(\lambda^\clubsuit, \lambda^\clubsuit, \lambda^\clubsuit)$,  $(\lambda^\spadesuit, \lambda^\spadesuit, \lambda^\spadesuit)$.

\vskip 1em
\noindent
$(\Leftarrow)$  As each member of $L(s,v)$ satisfies only a necessary condition for being incompressible (compare  (C) and Proposition \ref{irred}), it may happen that  $\lambda^\clubsuit$ or $\lambda^\spadesuit$ are compressible. By Proposition \ref{marozw},  the only possibility is that $\lambda^\spadesuit$ can be compressed to  $\lambda^\clubsuit$; that is, there is a homomorphism $h\colon \ka I(5)\to \ka I(4)$ such that 
$h_*(\lambda^\spadesuit)=\lambda^\clubsuit$. Consequently,  $|E(\lambda^\clubsuit)|\ge |E(\lambda^\spadesuit)|$. The latter  is impossible, as  $|E(\lambda^\clubsuit)|= 24$ while  $|E(\lambda^\spadesuit)|=26$. \hfill $\square$

\section{Isomorphic incompressible cliques.  Automorphisms}
From now on, we shall interpret $\text{Cq}^\clubsuit$ and $\text{Cq}^\spadesuit$ as families of 12-cliques rather than formal combinations of the 3-intervals these cliques consist of.  As in the preceding section, let $\aut s$ be  the automorphism group of  the graph $\ka I(s)$. Let $\Aut 3 s$ consists of the product mappings  $f=f_1\times f_2\times f_3$, where $f_i\in \aut s$. It is clear that for every $\ka C  \in \text{Cq}^\clubsuit$ and $f\in \Aut 3 4$, we have $f(\ka C)\in  \text{Cq}^\clubsuit$.  Therefore,  the group $\Aut 3 4$ acts  on $\text{Cq}^\clubsuit$. It can be rather easily computed that  $\text{Cq}^\clubsuit$ is an orbit of $\Aut 3 4$. Consequently,  we have:
\begin{pr}
All cliques belonging to $\operatorname{Cq^\clubsuit}$ are isomorphic. 
\end{pr}

Similarly, $\Aut 3 5$ acts on $\text{Cq}^\spadesuit$.  It appears that  $\text{Cq}^\spadesuit$  splits into $3$ orbits under the action of  $\Aut 3 5$. (Two of them are of cardinality 64 while the remaining orbit is of cardinality  128).  Clearly, it does not necessarily mean that there are three pairwise non-isomorphic cliques. In fact, they are not.  

Let $S_3$ be the symmetry group of $\{1,2,3\}$. Each $\sigma\in S_3$ induces the isometry $s_\sigma$ of $\er^3$ which in turn extends to intervals:  $s_{\sigma}(I)=I_{\sigma^{-1}(1)}\times I_{\sigma^{-1}(2)}\times I_{\sigma^{-1}(3)}$.  Obviously, $s_\sigma$  defines an isomorphism of cliques.  It is rather obvious that the composites $s_\sigma \circ f$, where $\sigma \in S_3$  and $f\in \Aut 3 s$,  form a group of isomorphisms of cliques.  Algebraically, it is a semidirect product of $S_3$ and $\Aut 3 s$.  Let us denote it by $A^3(s)$.  A computation shows  that $\text{Cq}^\spadesuit$ splits into  two orbits under the action  of $A^3(5)$. Let us fix some representatives of these orbits for further discussion, and call them $\ka D^1$ and $\ka D^2$ (Table \ref{representatives}, Figure \ref{representatives}).
\begin{table}[ht]
$$
\begin{array}{rccc}
1.\,\, & [0,1]\times[1,3]\times[2,3] & \,\,\,&[0, 1]\times[1,3]\times[2,3]\\
2.\,\, & [1,2]\times[2,3]\times[1,3] & \,\,\,&[1,2]\times[2,3]\times[1,3]\\ 
3.\,\, & [1,3]\times[1,2]\times[2,3] & \,\,\,&[1,3]\times[1,2]\times[2,3]\\ 
4.\,\, & [1,3]\times[2,3]\times[0,1] & \,\,\,&[1,3]\times[2,3]\times[0,1]\\
5.\,\, & [2,3]\times[0,1]\times[1,3] & \,\,\,&[2,3]\times[0,1]\times[1,3]\\ 
6.\,\, & [2,3]\times[1,3]\times[1,2] & \,\,\,&[2,3]\times[1,3]\times[1,2]\\ 
7.\,\, & [3,4]\times[3,5]\times[5,6] & \,\,\,&[3,4]\times[3,5]\times[4,5]\\
8.\,\, & [3,4]\times[4,5]\times[3,5] & \,\,\,&[3,4]\times[5,6]\times[3,5]\\
9.\,\, & [3,5]\times[5,6]\times[3,4] & \,\,\,&[3,5]\times[3,4]\times[5,6]\\
10.\,\,& [3,5]\times[3,4]\times[4,5]& \,\,\,&[3,5]\times[4,5]\times[3,4]\\
11.\,\,& [4,5]\times[3,5]\times[3,4]& \,\,\,&[4,5]\times[3,4]\times[3,5]\\
12.\,\,& [5,6]\times[3,4]\times[3,5]& \,\,\,&[5,6]\times[3,5]\times[3,4]
\end{array}.
$$
\caption{A system of representatives for $\text{Cq}^\spadesuit/A^3(5)$. \protect\label{representatives}} 
\end{table}

\begin{figure}[!ht]
\label{representativesf}
\centering
\includegraphics[width=0.4\textwidth]{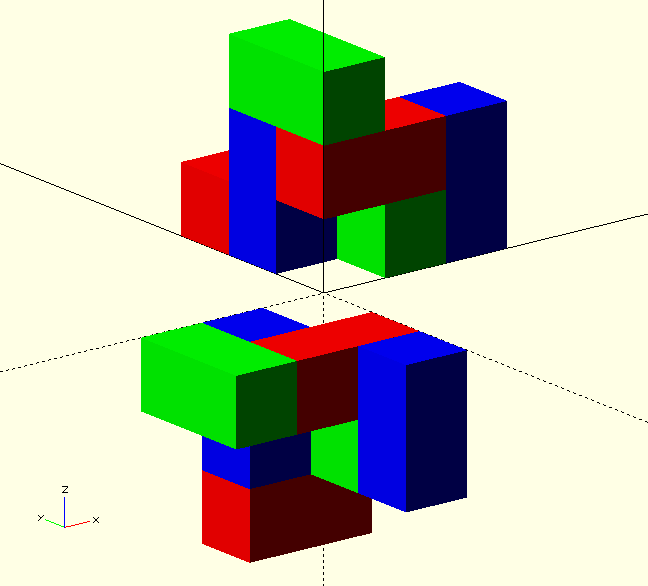}
\qquad
\includegraphics[width=0.4\textwidth]{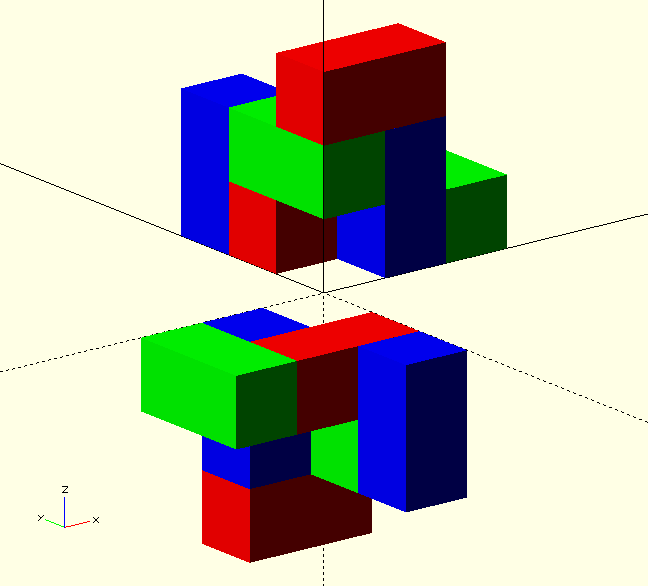}
\caption{Illustrations of cliques $\ka D^1$ and $\ka D^2$. They are translated so that their circumscribed boxes  are centred at the origin.}
\end{figure}

It seems to be a proper place to  give a formal definition of an isomorphism of cliques. 

For every pair of $n$-intervals $I$ and $J$ we define a $0/1$ vector $\varepsilon =\varepsilon(I,J)$  in $\er^n$ as follows
$$
\varepsilon_i=\left\{
\begin{array}{ll}
1, & \text{if $I_i$ and $J_i$ have exactly one point in common,} \\
0, & \text{otherwise.}
\end{array}
\right.  
$$
Let  $\ka C$ and $\ka D$ be two subfamilies of $\ka I^n$.  A bijection $f\colon \ka C\to \ka D$  is an \textit{isomorphism} between these families if  there is a permutation $\sigma \in S_n$ so that 
$$
s_\sigma(\varepsilon( f(I),  f(J)))=\varepsilon( s_\sigma(f(I)),  s_\sigma( f(J)))=\varepsilon(I, J)
$$ 
for every $I, J\in \ka C$. Obviously, if one of the two families is a clique, then the other is a clique as well.   Isomorphic  families will also be called \textit{combinatorially equivalent}.  

Observe that   if $\ka C\in  \text{Cq}^\clubsuit$ and $\ka D\in \text{Cq}^\spadesuit$, then they cannot be combinatorially equivalent;  otherwise, $|E(\lambda^\clubsuit|= |E(\lambda^\spadesuit)|$, which is not true.  Now, we are about to show  that cliques $\ka D^1$ and $\ka D^2$ are not equivalent as well. To this end,  it suffices to show that their automorphism groups are non-isomorphic. 

For every finite family $\ka D=\{I^1, \ldots, I^m\}\subset \ka I^n$, we may define its \textit{adjacency matrix} $A_{\ka D} =[\varepsilon_{ij}]$ so that 
$\varepsilon_{ij}= \varepsilon(I^i, I^j)$.  Clearly, the automorphism group $\aut{\ka D}$ is isomorphic to the group of all  these permutations $\pi \in S_m$ for whose there are $\sigma\in S_n$ such that $s_\sigma(\varepsilon_{\pi(i)\pi(j)})=\varepsilon_{ij}$ for every $i$ and $j$.  Table \ref{sasiedztwa} contains the adjacency matrices of $\ka D^1$ and $\ka D^2$; in order to simplify the notation  each $0/1$ wector $(\alpha, \beta, \gamma)$ is replaced by the number $\alpha  +2\beta + 4\gamma$.       

\begin{table}
$$
\left[
\begin{array}{cccccccccccc}
0& 1& 1& 1& 2& 4& 2& 4& 2& 4& 6& 6\\
1& 0& 2& 4& 1& 1& 2& 4& 2& 4& 6& 6\\
1& 2& 0& 2& 2& 4& 1& 5& 1& 5& 4& 4\\
1& 4& 2& 0& 4& 4& 3& 1& 3& 1& 2& 2\\
2& 1& 2& 4& 0& 2& 1& 5& 1& 5& 4& 4\\
4& 1& 4& 4& 2& 0& 3& 1& 3& 1& 2& 2\\
2& 2& 1& 3& 1& 3& 0& 4& 4& 2& 1& 4\\
4& 4& 5& 1& 5& 1& 4& 0& 2& 2& 1& 2\\
2& 2& 1& 3& 1& 3& 4& 2& 0& 4& 4& 1\\
4& 4& 5& 1& 5& 1& 2& 2& 4& 0& 2& 1\\
6& 6& 4& 2& 4& 2& 1& 1& 4& 2& 0& 1\\
6& 6& 4& 2& 4& 2& 4& 2& 1& 1& 1& 0
\end{array}
\right],
\qquad
\left[
\begin{array}{cccccccccccc}                                    
0& 1& 1& 1& 2& 4& 2& 4& 2& 4& 6& 6\\
1& 0& 2& 4& 1& 1& 2& 4& 2& 4& 6& 6\\
1& 2& 0& 2& 2& 4& 1& 5& 1& 5& 4& 4\\
1& 4& 2& 0& 4& 4& 3& 1& 3& 1& 2& 2\\
2& 1& 2& 4& 0& 2& 1& 5& 1& 5& 4& 4\\
4& 1& 4& 4& 2& 0& 3& 1& 3& 1& 2& 2\\
2& 2& 1& 3& 1& 3& 0& 2& 4& 4& 1& 4\\
4& 4& 5& 1& 5& 1& 2& 0& 4& 2& 1& 2\\
2& 2& 1& 3& 1& 3& 4& 4& 0& 2& 4& 1\\
4& 4& 5& 1& 5& 1& 4& 2& 2& 0& 2& 1\\
6& 6& 4& 2& 4& 2& 1& 1& 4& 2& 0& 1\\
6& 6& 4& 2& 4& 2& 4& 2& 1& 1& 1& 0
\end{array}
\right]
$$
\caption{Adjacency matrices of $\ka D^1$ and $\ka D^2$.  \protect\label{sasiedztwa}} 
\end{table}

The following simple idea can be applied in order to compute the automorphism group of  $\ka D$  efficiently.  Let $b_{ij}$ be the cardinality  of the set $\{k\colon \varepsilon_{ik}=\varepsilon_{jk}\}$.  Let $B_{\ka D}=[b_{ij}]$.  A permutation $\pi\in S_m$ is a \textit{protoautomorphism} of  $\ka D$ if it fixes $B_{\ka D}$; that is, $[b_{ij}]=[b_{\pi(i)\pi(j)}]$. Clearly, each automorphism is a protoautomorphism.  On the other hand, it is much easier to determine whether a permutation is a protoautomorphism than an automorphism.  It appears that both $\ka D^1$ and $\ka D^2$ have  only 48 protoautomorphism.  Finally, it remains to single out automorphisms from protoautomorphisms, which is an easy task as the number of protoautomorphisms is small. Again, it appears that both automorphism groups are of the same order 24.  All automorphisms of both cliques are collected in Tabel \ref{auto}.  
\begin{table}
{\footnotesize
$$
\begin{array}{l}
(1)\, (2)\, (3)\, (4)\, (5)\, (6)\, (7)\, (9)\, (10)\, (11)\, (12)\, (8)\\
(1)\, (2)\, (3)\, (4)\, (5)\, (6)\, (7, 9)\, (10, 8)\, (11, 12)\\
(1, 2)\, (3, 5)\, (4, 6)\, (7)\, (9)\, (10)\, (11)\, (12)\, (8)\\
(1, 2)\, (3, 5)\, (4, 6)\, (7, 9)\, (10, 8)\, (11, 12)\\
(1, 3, 4, 2, 5, 6)\, (7, 11, 10, 9, 12, 8)\\
(1, 3, 4, 2, 5, 6)\, (7, 12, 10)\, (9, 11, 8)\\
(1, 4, 5)\, (2, 6, 3)\, (7, 8, 12, 9, 10, 11)\\
(1, 4, 5)\, (2, 6, 3)\, (7, 10, 12)\, (9, 8, 11)\\
(1, 5, 4)\, (2, 3, 6)\, (7, 11, 10, 9, 12, 8)\\
(1, 5, 4)\, (2, 3, 6)\, (7, 12, 10)\, (9, 11, 8)\\
(1, 6, 5, 2, 4, 3)\, (7, 8, 12, 9, 10, 11)\\
(1, 6, 5, 2, 4, 3)\, (7, 10, 12)\, (9, 8, 11)\\
(1, 7)\, (2, 9)\, (3, 8)\, (4, 12)\, (10, 5)\, (11, 6)\\
(1, 7, 2, 9)\, (3, 8, 5, 10)\, (12, 6, 11, 4)\\
(1, 8, 2, 10)\, (3, 12, 5, 11)\, (9, 6, 7, 4)\\
(1, 8)\, (2, 10)\, (3, 12)\, (4, 9)\, (11, 5)\, (6, 7)\\
(1, 9, 2, 7)\, (3, 10, 5, 8)\, (12, 4, 11, 6)\\
(1, 9)\, (2, 7)\, (3, 10)\, (4, 11)\, (12, 6)\, (5, 8)\\
(1, 10)\, (2, 8)\, (3, 11)\, (4, 7)\, (9, 6)\, (12, 5)\\
(1, 10, 2, 8)\, (3, 11, 5, 12)\, (9, 4, 7, 6)\\
(1, 11)\, (2, 12)\, (3, 7)\, (4, 8)\, (9, 5)\, (10, 6)\\
(1, 11, 2, 12)\, (3, 7, 5, 9)\, (10, 4, 8, 6)\\
(1, 12, 2, 11)\, (3, 9, 5, 7)\, (10, 6, 8, 4)\\
(1, 12)\, (2, 11)\, (3, 9)\, (4, 10)\, (5, 7)\, (6, 8)\\
\end{array}
\qquad
\begin{array}{l}
(1)\, (2)\, (3)\, (4)\, (5)\, (6)\, (7)\, (9)\, (10)\, (11)\, (12)\, (8)\\
(1)\, (2)\, (3)\, (4)\, (5)\, (6)\, (7, 9)\, (10, 8)\, (11, 12)\\
(1, 2)\, (3, 5)\, (4, 6)\, (7)\, (9)\, (10)\, (11)\, (12)\, (8)\\
(1, 2)\, (3, 5)\, (4, 6)\, (7, 9)\, (10, 8)\, (11, 12)\\
(1, 3, 4, 2, 5, 6)\, (7, 11, 10)\, (9, 12, 8)\\
(1, 3, 4, 2, 5, 6)\, (7, 12, 10, 9, 11, 8)\\
(1, 4, 5)\, (2, 6, 3)\, (7, 8, 11, 9, 10, 12)\\
(1, 4, 5)\, (2, 6, 3)\, (7, 10, 11)\, (9, 8, 12)\\
(1, 5, 4)\, (2, 3, 6)\, (7, 11, 10)\, (9, 12, 8)\\
(1, 5, 4)\, (2, 3, 6)\, (7, 12, 10, 9, 11, 8)\\
(1, 6, 5, 2, 4, 3)\, (7, 8, 11, 9, 10, 12)\\
(1, 6, 5, 2, 4, 3)\, (7, 10, 11)\, (9, 8, 12)\\
(1, 7, 3, 12, 4, 10, 2, 9, 5, 11, 6, 8)\\
(1, 7, 5, 11, 4, 10)\, (2, 9, 3, 12, 6, 8)\\
(1, 8, 6, 11, 5, 9, 2, 10, 4, 12, 3, 7)\\
(1, 8, 4, 12, 5, 9)\, (2, 10, 6, 11, 3, 7)\\
(1, 9, 5, 12, 4, 8)\, (2, 7, 3, 11, 6, 10)\\
(1, 9, 3, 11, 4, 8, 2, 7, 5, 12, 6, 10)\\
(1, 10, 4, 11, 5, 7)\, (2, 8, 6, 12, 3, 9)\\
(1, 10, 6, 12, 5, 7, 2, 8, 4, 11, 3, 9)\\
(1, 11)\, (2, 12)\, (3, 8)\, (4, 7)\, (9, 6)\, (10, 5)\\
(1, 11, 2, 12)\, (3, 8, 5, 10)\, (9, 4, 7, 6)\\
(1, 12, 2, 11)\, (3, 10, 5, 8)\, (9, 6, 7, 4)\\
(1, 12)\, (2, 11)\, (3, 10)\, (4, 9)\, (5, 8)\, (6, 7)\\
\end{array}
$$
}
\caption{ The isomorphic copies of $\aut {\ka D^1}$ and  $\aut{\ka D^2}. $\protect\label{auto}}
\end{table}
Clearly,  groups $\aut{\ka D^1}$ and $\aut{\ka D^2}$ are different, as only the latter contains elements of order 12. In fact, we can easily identify these groups using GAP.  The first of them is an (external) semidirect product of cyclic group $\zet_3$ by dihedral group $\operatorname{Dih}_4$ of order 8, where the action of the latter on the former is given by a homomorphism  whose kernel is one of the Klein four-subgroups  of the dihedral group.  The second is the simple product  of $\zet_3$ and $\operatorname{Dih}_4$.  
\begin{pr}
The quotient space $\operatorname{Cq}^\spadesuit\!/\!\!\approx$, where $\approx$ is the combinatorial equivalence of cliques, consists of two classes
\end{pr}  

Since all cliques belonging to $\operatorname{Cq}^\clubsuit$ are isomorphic, their automorphism groups are isomorphic as well. Therefore, it suffices to fix any $\ka C\in \operatorname{Cq}^\clubsuit$, and find $\aut{\ka C}$ in order to know the structure of the automorphism groups of the remaining cliques. We could take for example  clique $\ka C$ described in  Example \ref{np12klika}. By much the same argument as in the case of $\ka D^i$ one finds that $\aut{\ka C}$ is of order 48.  (Interestingly enough, the protoautomorphism group has in this case 3070 elements).  The following permutations are generators of the isomorphic copy of $\aut{\ka C}$: $(1, 12, 7, 4, 9, 6)\, (2, 10, 5, 3, 11, 8)$, $(1, 12, 4, 9)\, (2, 11, 3, 10)\, (5, 6, 8, 7)$.  Since $|A^3(s)|=|S_3||\aut s|^3= 6 \cdot 8^3$,  for $s\ge 4$,   $|\aut{\ka C}|=48$ and $|\operatorname{Cq}^\clubsuit|=64$, it follows that each automorphism of $\ka C$ is the restriction to $\ka C$ of an element from $A^3(4)$. Similarly, all automorphisms of  $\ka D^1$ and $\ka D^2$  can be identified with corresponding elements of $A^3(5)$. 

\section{Compressible cliques.}
We have classified all incompressible 12-cliques in $\ka I^3$ up to combinatorial equivalence, however, it is not the whole picture.  There are 12-cliques which are of different combinatorial type from the three described so far. Clearly, if $\ka C$ is such a clique, then  there is a homomorphism 
$f=f_1\times f_2\times f_2\colon \ka I^3\to \ka I(s)^3$, where $s\in \{4,5\}$ such that  ${\ka D}=f(\ka {C})$ is an incompressible clique  (see Section \ref{mincl}). Since $\varepsilon(I,J)\le \varepsilon(f(I),f(J))$,  for every $I, J\in \ka I^3$, we deduce that   $A_{\ka C}\le A_{\ka D}$, where both inequalities are stated with respect to the coordinatewise order.  Since we have only three combinatorial types of incompressible cliques, we deduce that the number of combinatorial types of 12-cliques is  finite.  As we briefly explain in this section, this number is 5.  Details will be published elsewhere.

It can be shown that if $A_{\ka C}\neq A_{\ka D}$, then $\ka D \in \operatorname{Cq}^\spadesuit$.   As our objective is to characterize all 12-cliques up to combinatorial equivalence, we may assume that $\ka D\in \{\ka D^1, \ka D^2\}$. Let us assume that the intervals belonging to $\ka D=\ka D^i$ are labelled as in Table \ref{representatives}. Moreover, let the elements of $\ka C=\ka C^i$ be labelled so that $I=f(J)\in \ka D$ and $J\in \ka C$ have the same number. If we take into account authomorphisms of $\ka D$, the fact that  $\varepsilon(I^1, I^8)=(1,0,1)$ and the assumption $A_{\ka C} \neq A_{\ka D}$, then we can set that $\varepsilon(J^1, J^8)=(0,0,1)$.  It appears that this equality determines $A_\ka C$ for both $\ka D^1$ and $\ka D^2$.  In Table \ref{sasiedztwa2}, the  adjacency matrices of $\ka C^1$ and $\ka C^2$ encoded in the same manner as in Table \ref{sasiedztwa} are presented. The subsequent table shows the corresponding 12-cliques.

\begin{table}
$$
\left[
\begin{array}{cccccccccccc}
0& 1& 1& 1& 2& 4& 2& 4& 2& 4& 6& 6\\
1& 0& 2& 4& 1& 1& 2& 4& 2& 4& 6& 6\\
1& 2& 0& 2& 2& 4& 1& 4& 1& 4 & 4& 4\\
1& 4& 2& 0& 4& 4& 2& 1& 2& 1& 2& 2\\
2& 1& 2& 4& 0& 2& 1& 4& 1& 4& 4& 4\\
4& 1& 4& 4& 2& 0& 2& 1& 2& 1& 2& 2\\
2& 2& 1& 2& 1& 2& 0& 4& 4& 2& 1& 4\\
4& 4& 4& 1& 4& 1& 4& 0& 2& 2& 1& 2\\
2& 2& 1& 2& 1& 2& 4& 2& 0& 4& 4& 1\\
4& 4& 4& 1& 4& 1& 2& 2& 4& 0& 2& 1\\
6& 6& 4& 2& 4& 2& 1& 1& 4& 2& 0& 1\\
6& 6& 4& 2& 4& 2& 4& 2& 1& 1& 1& 0
\end{array}
\right],
\qquad
\left[
\begin{array}{cccccccccccc}                                    
0& 1& 1& 1& 2& 4& 2& 4& 2& 4& 6& 6\\
1& 0& 2& 4& 1& 1& 2& 4& 2& 4& 6& 6\\
1& 2& 0& 2& 2& 4& 1& 4& 1& 4& 4& 4\\
1& 4& 2& 0& 4& 4& 2& 1& 2& 1& 2& 2\\
2& 1& 2& 4& 0& 2& 1& 4& 1& 4& 4& 4\\
4& 1& 4& 4& 2& 0& 2& 1& 2& 1& 2& 2\\
2& 2& 1& 2& 1& 2& 0& 2& 4& 4& 1& 4\\
4& 4& 4& 1& 4& 1& 2& 0& 4& 2& 1& 2\\
2& 2& 1& 2& 1& 2& 4& 4& 0& 2& 4& 1\\
4& 4& 4& 1& 4& 1& 4& 2& 2& 0& 2& 1\\
6& 6& 4& 2& 4& 2& 1& 1& 4& 2& 0& 1\\
6& 6& 4& 2& 4& 2& 4& 2& 1& 1& 1& 0
\end{array}
\right]
$$
\caption{Adjacency matrices of $\ka C^1$ and $\ka C^2$.  \protect\label{sasiedztwa2}} 
\end{table}

\begin{table}
$$
\begin{array}{rccc}
1.\,\, & [0,1]\times[1,3]\times[2,3] & \,\,\,&[0, 1]\times[1,3]\times[2,3]\\
2.\,\, & [1,2]\times[2,3]\times[1,3] & \,\,\,&[1,2]\times[2,3]\times[1,3]\\ 
3.\,\, & [1, \frac 5 2]\times[1,2]\times[2,3] & \,\,\,&[1,\frac 5 2]\times[1,2]\times[2,3]\\ 
4.\,\, & [1,3]\times[2,3]\times[0,1] & \,\,\,&[1,3]\times[2,3]\times[0,1]\\
5.\,\, & [2, \frac 5 2]\times[0,1]\times[1,3] & \,\,\,&[2,\frac 5 2]\times[0,1]\times[1,3]\\ 
6.\,\, & [2,3]\times[1,3]\times[1,2] & \,\,\,&[2,3]\times[1,3]\times[1,2]\\ 
7.\,\, & [\frac 5 2, 4]\times[3,5]\times[5,6] & \,\,\,&[\frac 5 2, 4]\times[3,5]\times[4,5]\\
8.\,\, & [3,4]\times[4,5]\times[3,5] & \,\,\,&[3,4]\times[5,6]\times[3,5]\\
9.\,\, & [\frac 5 2, 5]\times[5,6]\times[3,4] & \,\,\,&[\frac 5 2 ,5]\times[3,4]\times[5,6]\\
10.\,\,& [3,5]\times[3,4]\times[4,5]& \,\,\,&[3,5]\times[4,5]\times[3,4]\\
11.\,\,& [4,5]\times[3,5]\times[3,4]& \,\,\,&[4,5]\times[3,4]\times[3,5]\\
12.\,\,& [5,6]\times[3,4]\times[3,5]& \,\,\,&[5,6]\times[3,5]\times[3,4]
\end{array}.
$$
\caption{Compressible cliques $\ka C^1$, $\ka C^2$. \protect\label{repr2}}  
\end{table}
It can be shown that $\ka C^1$ and $\ka C^2$ are not isomorphic, despite the fact that their automorphism groups are both isomorphic to $\operatorname{Dih}_4$.

\begin{figure}[!ht]
\centering
\includegraphics[width=0.45\textwidth]{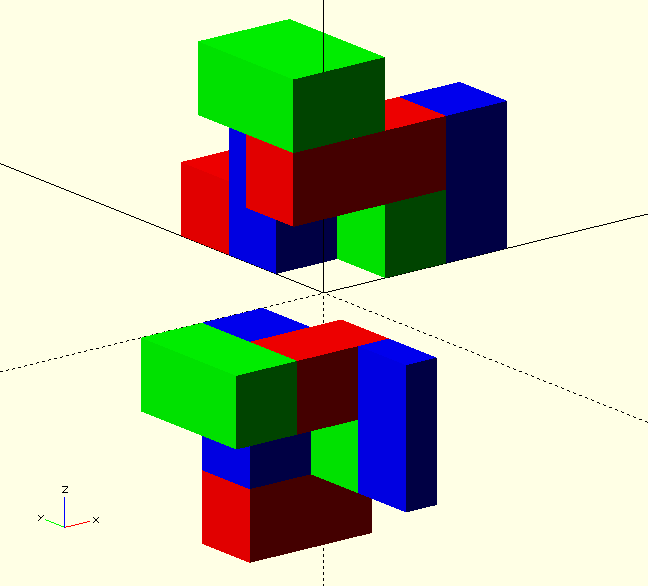}
\qquad
\includegraphics[width=0.45\textwidth]{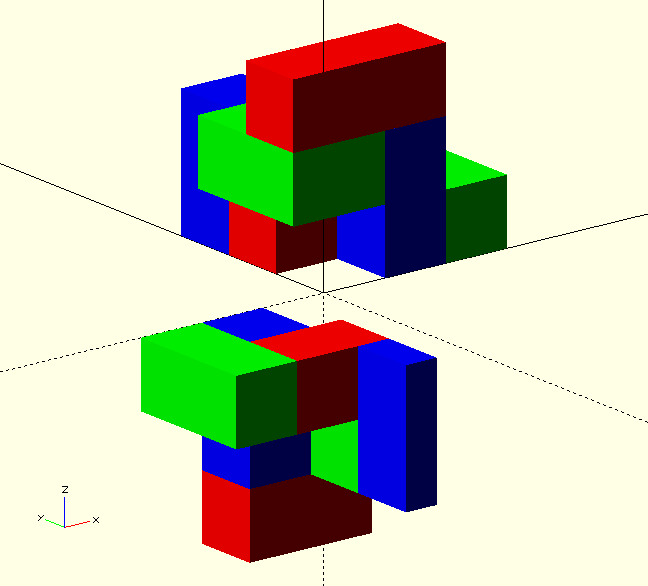}
\caption{Illustrations of compressible cliques $\ka C^1$ and $\ka C^2$. They are translated so that their circumscribed boxes  are centred at the origin.}
\end{figure}

\section{Isometric incompressible cliques.  Chirality}

Let $I(s)=[0, s+1]\times [0,s+1]\times [0,s+1]$.  Let $\operatorname{Iso}^3(s)$ be the group of isometries of the cube $I(s)$.  For  $s=4$, $\operatorname{Iso}^3(s)$ acts on $\operatorname{Cq}^\clubsuit$, while for $s=5$, $\operatorname{Iso}^3(s)$ acts on  $\operatorname{Cq}^\spadesuit$. The action is defined naturally: If $U\in \operatorname{Iso}^3(s)$ and $\ka C$ is an incompressible 12-clique, then $U(\ka C)=\{U(I) \colon I \in \ka C\}$.  Let us remark, that $\operatorname{Iso}^3(s)$ can be identified with the subgroup  of $A^3(s)$ consisting of all the  elements $s_\sigma\circ f=s_\sigma\circ (f_1\times f_2\times f_3)$ such that $\sigma$ is an arbitrary permutation belonging to $S_3$ while   each $f_i$  is either the identity mapping or the reflection about $(s+1)/2$; that is, $f_i([a,b])= [s+1-b, s+1-a]$, for $[a,b ]\in \ka I(s)$. 

Each group $\operatorname{Iso}^3(s)$ contains the subgroup $\operatorname{Iso}^3_+(s)$ consisting of  preserving-orientation isometries.    
Since incompressible cliques of intervals are geometric objects, it makes sense to characterize them up to congruency; that is, to find a system of representatives for each of the quotient spaces: $\operatorname{Cq}^\clubsuit/\operatorname{Iso}^3_+(4)$,  $\operatorname{Cq}^\clubsuit/\operatorname{Iso}^3(4)$,  $\operatorname{Cq}^\spadesuit/\operatorname{Iso}^3_+(5)$,  $\operatorname{Cq}^\spadesuit/\operatorname{Iso}^3(5)$.  It appears that some of the cliques are chiral while the other are achiral.  In our context, a clique $\ka C$ is \textit{achiral} if its orbits with respect to  $\operatorname{Iso}^3_+(s)$ and $\operatorname{Iso}^3(s)$ coincide; otherwise, it is  \textit{chiral}.   

\begin{figure}
\label{labeling}
\centering
\includegraphics{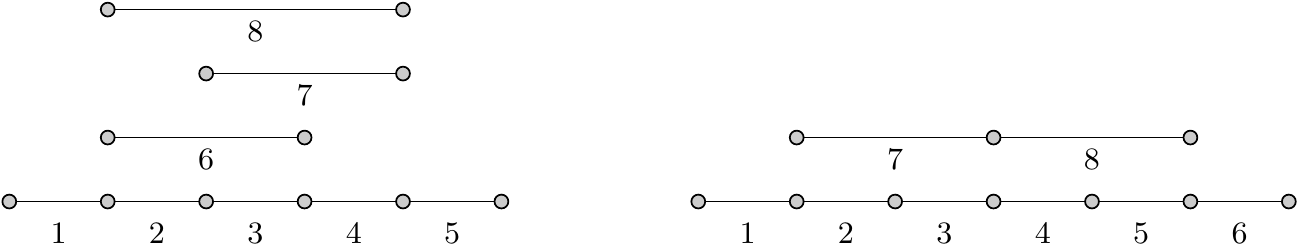}
\caption{Labelings of $\supp\lambda^\clubsuit$ and $\supp\lambda^\spadesuit$.}
\end{figure}

 To simplify the exposition, we label the elements of $\supp \lambda^\clubsuit$ as follows:  $[i-1,i]\mapsto i$, for $i=1,\ldots, 5$; $[1,3]\mapsto 6$, $[2,4]\mapsto 7$, $[1,4]\mapsto 8$.  Then we extend this labelling to 3-intervals; for example,  interval   $I=[0,1]\times[3,4]\times [1,4]$ is labeled by $148$.  Similarly we label the elements of of $\supp \lambda^\spadesuit$:
$[i-1,i] \mapsto i$, for $i=1,\ldots, 6$; $[1,3]\mapsto 7$, $[3,5]\mapsto 8$ (Figure \ref{labeling}).  We gather the information concerning the quotient spaces under the discussion in Tables \ref{orbity_4}--\ref{orbity_5}

\begin{table}[!ht]
\begin{tabular}{|ccc|}
\hline
No. & orbit representative & orbit length \\
\hline
1.                 
& [148, 157, 246, 253, 624, 814, 325, 715, 462, 481, 532, 571]  & 8\\ 
2.                 
& [165, 134, 285, 274, 651, 852, 341, 742, 413, 427, 516, 528]	 & 8\\
\hline

3.                 &     [116, 128, 213, 227, 652, 842, 351, 741, 465, 434, 585, 574]  &  24\\
4.                  &    [116, 123, 218, 227, 645, 855, 344, 754, 482, 471, 562, 531]  &  24\\
\hline
\end{tabular}

\caption{A system of representatives for $\operatorname{Cq}^\clubsuit/\operatorname{Iso}^3_+(4)$. Two blocks of cliques belonging to the same orbit of $\operatorname{Iso}^3(4)$ are distinguished. All cliques are chiral. \protect\label{orbity_4}}
\end{table}

\begin{table} [!ht]
\begin{tabular}{|ccc|}
\hline
No. & orbit representative & orbit length \\
\hline
 I.  &  class of  the combinatorial equivalence $\approx$ &\\
\hline
1.   & [177, 233, 727, 732, 313, 371, 445, 464, 886, 858, 548, 684]  & 12\\
2.   & [177, 233, 723, 772, 317, 331, 446, 454, 885, 868, 584, 648]  & 12\\
\hline
3.   & [177, 233, 727, 732, 313, 371, 446, 468, 885, 854, 544, 688]  & 12\\
4.   & [177, 233, 723, 772, 317, 331, 486, 464, 845, 858, 544, 688]  & 12\\
\hline
5.   &[177, 233, 727, 732, 313, 371, 485, 454, 846, 868, 588, 644]   & 12\\
6.   &[177, 233, 723, 772, 317, 331, 445, 458, 886, 864, 588, 644]   & 12\\
\hline
7.   & [177, 233, 727, 732, 313, 371, 486, 458, 845, 864, 584, 648]  & 24\\
8.   & [177, 233, 723, 772, 317, 331, 485, 468, 846, 854, 548, 684]  & 24\\
\hline
9.   & [173, 237, 723, 731, 317, 372, 486, 458, 845, 864, 584, 648]  &  4\\
10. & [137, 273, 713, 732, 327, 371, 485, 468, 846, 854, 548, 684]  &  4\\ 
\hline
 II.  &  class of  the combinatorial equivalence $\approx$ &\\
\hline
11. & [177, 233, 727, 732, 313, 371, 445, 458, 886, 864, 588, 644]  & 24\\
12. & [177, 233, 727, 732, 313, 371, 446, 454, 885, 868, 584, 648]  & 24\\
\hline
13. & [177, 233, 727, 732, 313, 371, 485, 468, 846, 854, 548, 684]  &  24\\
14. & [177, 233, 723, 772, 317, 331, 486, 458, 845, 864, 584, 648]  & 24\\
\hline
15. & [177, 233, 727, 732, 313, 371, 486, 464, 845, 858, 544, 688]  & 24 \\
\hline
16. & [173, 237, 723, 731, 317, 372, 485, 468, 846, 854, 548, 684]  & 8   \\
\hline
\end{tabular}
\caption{A system of representatives for $\operatorname{Cq}^\spadesuit/\operatorname{Iso}^3_+(5)$. There are two blocks of cliques corresponding to two classes of 
the combinatorial equivalence. Each block consists of  subblocks. Cliques within a subblock belong to the same orbit of   $\operatorname{Iso}^3(5)$. Cliques  15 and 16 are achiral.\protect\label{orbity_5}}
\end{table}

\newpage

\section{Appendix A: $\aut s$}

\begin{figure}[htb]
\centering
\includegraphics{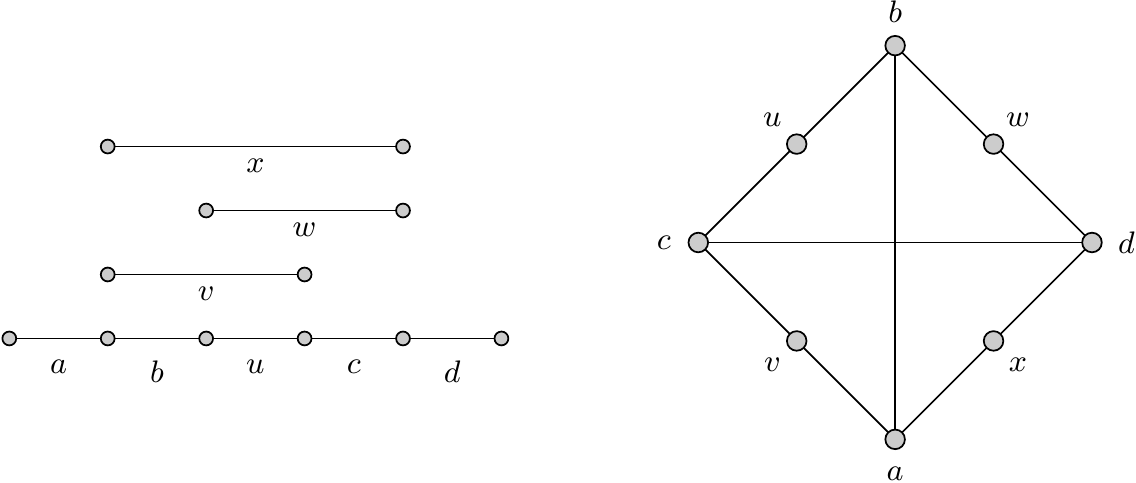}
\caption{ Diagrams for $\ka I(4)$.\protect{\label{s4}}} 
\end{figure}

Our goal  is to describe the group $\aut s$ of all automorphisms of the graph $\ka I (s)$ (see: Section \ref{upper bound} for the definition of $\ka I(s)$).
If $s\in \{1,2\}$, then $\ka I (s)$ is a path of length $s$  and consequently  $\aut s $ is isomorphic to $\zet_2$.  If $s=3$, then $\ka I(s)$ is a cycle of length $5$. Therefore, $\aut 3$ is isomorphic to $\operatorname{Dih}_5$, which can be identified with the group of isometries of a regular pentagon.  In the case  of $s=4$ ,  one can see that $\aut s$ is isomorphic to $\operatorname{Dih_4}$,  the group of isometries of a square.  Figure \ref{s4} provides a sufficient explanation.
 
\noindent       
In fact, the latter statement remains valid for arbitrary $s\ge 4$:
\begin{tw}
$\aut s \cong \operatorname{Dih}_4$, for every  $s\ge 4$.
\end{tw}
 We are going to prove this result straightforwardly; that is,  by determining all the elements of $\aut s$.  We begin with analysing the degrees of the vertices of  $\ka I(s)$.  We adopt the following shorthand notation:  Each vertex (interval) $[i,j]$  of the graph $\ka I(s)$ is denoted by $i\dkrop j$. 
 \begin{lemat}
Let $s\ge 4$. Let $i\dkrop j \in V(\ka I(s))$. Then  $\deg (i\dkrop j) $ is equal to:
\begin{enumerate}
\item
$s-(j-i)-1$,  if $2\le i$ and  $j\le s-1$;
\item
$s-(j-i)$, if  $i$, $j$ satisfy one of the following relations:   $i=0$ and $ j=1$, $i=1$ and  $j\le s-1$, $2\le i$ and  $j=s$,  $i=s$ and $j=s+1$ ;
\item
$s-(j-i)+1=2$, if $i=1$ and $j=s$. 
\end{enumerate}
\end{lemat}
The proof is obvious,  the degree $\deg(i\dkrop j)$ is equal to the number of the intervals in  $\ka I(s)$ adjacent to $i\dkrop j$.  For example, if $i=1$ and  $j\le s-1$, then the only interval adjacent to $i\dkrop j$ from the left is $0\dkrop 1$, while  the intervals $j\dkrop k$, where $j+1\le k\le s$, are  adjacent from the right.  Therefore,  there are $1+(s-j)= s-(j-i)$ adjacent intervals in total, a result which agrees with  our lemma.   

\Skok
\noindent
\textit{Proof of the theorem.} 
It follows from the lemma that if $s\ge 4$, then the only vertices of degree $s-1$ in $\ka I(s)$ are $a= 0\dkrop 1$, $b=1\dkrop 2$, $c=s-1\dkrop s$ and $d= s\dkrop s+1$.  Thus, every automorphism of $\ka I(s)$ has to send these vertices onto themselves. Since the subgraph $G$ of $\ka I(s)$  induced by  $V=\{a,b,c,d\}$ is a disjoint union of two paths of length 1,  its automorphism group $\aut G $ is isomorphic to $\operatorname{Dih}_4$. We shall prove that every automorphism of  $G$ extends to an automorphism of  $\ka I(s)$ in a unique way.  Let us begin with listing  the automorphisms of $G$ .  
$$
\begin{array}{c}
\begin{array}{c c} 
\begin{array}{|c||c|c|c|c|}\hline
1 & a & b & c & d\\\hline
2 & d & c & a & b\\\hline
3 & b & a & d & c\\\hline
4 & c & d & b & a\\\hline
5 & b & a & c & d\\\hline
6 & a & b & d & c\\\hline
7 & d & c & b & a\\\hline
8 & c & d & a & b\\\hline
\end{array} 
&\qquad
\begin{array}{|c||r@{\dkrop}l|r@{\dkrop}l|r@{\dkrop}l|r@{\dkrop}l|}\hline
1 & 0& 1 & 1& 2 & s-1& s & s & s+1\\\hline
2 & s& s+1 & s-1& s & 0& 1 &  1 & 2\\\hline
3 & 1& 2 & 0& 1 & s& s+1 & s-1 & s\\\hline
4 & s-1& s & s& s+1 & 1& 2 &  0  & 1\\\hline
5 & 1& 2 & 0& 1 & s-1& s &  s & s+1\\\hline
6 & 0& 1 & 1& 2 & s& s+1 & s-1 & s\\\hline
7 & s& s+1 & s-1& s & 1& 2 & 0 & 1\\\hline
8 & s-1& s & s& s+1 & 0& 1 & 1 & 2\\\hline
\end{array}\\
\end{array}\\
\\
\text{Lists of the automorphisms of $G$ in both notations for vertices.}\\
\end{array}
$$
These lists  encode automorphisms in a rather obvious manner; for example, the automorphism defined by the row no. 5 acts as follows: $a\mapsto b$, $b\mapsto a$, $c\mapsto c$, $d\mapsto d$.   
 
Let us consider the vertices of the lowest degree in $\ka I(s)$, $s\ge 4$. These are $u=2\dkrop s-1$, $v=1\dkrop s-1$, $w=2\dkrop s$ oraz $x=1\dkrop s$.  The graph $H$  induced by the vertices $a$, $b$, $c$, $d$, $u$, $v$, $w$, $x$ is isomorphic to $\ka  I(4)$ (for $s=4$ they simply coincide) and Figure \ref{s4} represents $H$ equally well as it represents $\ka I(4)$. 

Clearly, the automorphisms of $G$ extend to the automorphisms of $H$ in a natural and unique way.  We arrange the automorphisms of $H$ into two tables.  The second table describes the action of these automorphisms restricted to the vertices $u$, $v$, $w$, $x$ in the original notation.  
$$
\begin{array}{c} 
\begin{array}{cc}
\begin{array}{|c||c|c|c|c|c|c|c|c|}\hline
1 & a & b & c & d & u & v & w & x\\\hline
2 & d & c & a & b & v & x & u & w\\\hline
3 & b & a & d & c & x & w & v & u\\\hline
4 & c & d & b & a & w & u & x & v\\\hline
5 & b & a & c & d & v & u & x & w\\\hline
6 & a & b & d & c & w & x & u & v\\\hline
7 & d & c & b & a & u & w & v & x\\\hline
8 & c & d & a & b & x & v & w & u\\\hline
\end{array}
&\qquad
\begin{array}{|c||r@{\dkrop}l|r@{\dkrop}l|r@{\dkrop}l|r@{\dkrop}l|}\hline
1 & 2 & s-1 & 1 & s-1 & 2 & s & 1 & s\\\hline
2 & 1 & s-1 & 1 & s & 2 & s-1 & 2 & s\\\hline
3 & 1 & s & 2 & s & 1 & s-1 & 2 & s-1\\\hline
4 & 2 & s & 2 & s-1 & 1 & s & 1 & s-1\\\hline
5 & 1 & s-1 & 2 & s-1 & 1 & s & 2 & s\\\hline
6 & 2 & s & 1 & s& 2 & s-1 & 1 & s-1\\\hline
7 & 2 & s-1 & 2 & s & 1 & s-1 & 1 & s\\\hline
8 & 1 & s & 1 & s-1 & 2 & s & 2 & s-1\\\hline
\end{array} 
\end{array}
\\
\\
\text{$\aut H$}\\
\end{array}
$$ 
Now, we are prepared to describe  all the automorphisms of  $\ka I(s)$, $s\ge 4$.  Since, as we shall show, they are uniquely determined by the automorphisms of $G$, it suffices to describe their action on vertices belonging to $\ka I(s)\setminus V(H)$.  These are the vertices of degree greater than 2 and smaller than $s-1$.  
The indices $i$, $j$ appearing in our description have to satisfy the following constraints:
$$
2<i<s-1;\quad  \text{$i<j$,  if both indices appear};\quad  2< j<s-1.
$$ 
The enumeration of automorphisms is consistent with the preceding tables. 
$$
\begin{array}{c} 
\begin{array}{|c||r@{\dkrop}l|r@{\dkrop}l|r@{\dkrop}l|r@{\dkrop}l|r@{\dkrop}l|}\hline
1 & 1& i            & 2& i            &   j & s-1     &   j& s        & i& j \\\hline
2 & s+1-i& s     &  s+1-i& s-1  & 1& s+1-j   &  2& s+1-j  & s+1-j& s+1-i \\\hline
3 & 2& i           &  1& i            &  j& s          &   j& s-1     &    i& j \\\hline
4 & s+1-i& s-1 &  s+1-i& s      &  2& s+1-j  &  1& s+1-j  &   s+1-j& s+1-i\\\hline
5 & 2& i           &  1& i            &  j & s-1     &   j& s         &   i& j \\\hline
6 & 1& i           &  2& i            &   j& s        &  j& s-1       &   i& j \\\hline
7 & s+1-i& s     &  s+1-i& s-1  & 2& s+1-j  &  1& s+1-j   &   s+1-j& s+1-i\\\hline
8 & s+1-i& s-1  &  s+1-i& s     & 1& s+1-j  &  2& s+1-j   &   s+1-j& s+1-i\\\hline
\end{array}\\
\\
\text{The action of $\aut s$ on  $\ka I(s)\setminus V(H)$.}
 \\
\end{array}
$$ 
The fact that the described eight mappings, call them $\varphi_i$, $i=1,\ldots, 8$, are automorphisms  of the graph $\ka I(s)$, $s\ge 4$,  is easily verified. As they extend the automorphism  of $G$, it remains to be shown that no other extensions exist.  Let $\varphi$ be any of the automorphisms.  Then there is a unique $i$,  such that  $\varphi|V(G)=\varphi_i|V(G)$.  As we have already mentioned, the extension  of $\varphi_i|V(G)$ to $V(H)$ is unique; therefore, $\varphi|V(H)=\varphi_i|V(H)$.  Since $H=\ka I(s)$ for $s=4$, we may further assume $s\ge 5$.  Observe now, that for every $y\in V(G)$, and  $2<\delta<s-1$ there is exactly one neighbouring  vertex $z$ of  degree $\delta$.   Let $y'=\varphi_i(y)=\varphi(y) $.  Since $y'\in V(G) $, there is only one vertex $z'$ being a neighbour  of $y'$ whose degree is $\delta$. This implies that  $z'=  \varphi(z)=\varphi_i(z)$. Thus, $\varphi$ and $\varphi_i$ coincide on the neighbourhood of $V(G)$.    Now, if a vertex $t$ does not belong to the neighbourhood of  $V(G)$, then there are $i$, $j$ such that  $t=i\dkrop j$ and $2<i<j<s-1$. Take $z_1=1\dkrop i$ oraz $z_2=j\dkrop s$.  These are neighbours of $t$. Moreover, they belong to the neighbourhood of $V(G)$, as the pairs $z_1$, $a$ and  $z_2$, $d$ are adjacent.  As a result,   both automorphisms  send $z_1$ on the same element  $z_1'$ and  $z_2$ onto $z'_2$.  Since $z_1$, $z_2$ have a unique common neighbour $t$,  the elements  $z_1'$, $z_2'$ have a unique common member $t'$.  Clearly, $t'= \varphi(t)=\varphi_i(t)$, which  completes the proof that $\varphi$ and $\varphi_i$ coincide.   \hfill $\square$    

\section{Appendix B:  \texorpdfstring{$\hat{L}(s,v)$} { L(s,v)} }
We fulfill our declaration made in Section \ref{mincl} to collect all  $\hat{L}(s,v)$ for $v=12,13$.  We explain how to read these tables taking row no. 9 of Table \ref{5,12} as an example:  Entries of this row bounded by double lines encode $\gamma=\lambda^*$, the combination  discussed in Proposition \ref{bezrozw}.  The fact that $e=29'$ means that the graph $G_\gamma$ has $29$ edges and is of type II.    
\begin{table}
$$
\begin{array}{|c||c|c|c|c|c||c|}\hline
\text{No.} & 0\dkrop 1& 1\dkrop 2&1\dkrop 3&2\dkrop 3&3\dkrop 4& e\\\hline
1&2 & 2 & 2 & 3 & 3&29\\\hline
\end{array} 
 $$ 
\caption{ $\hat{L}(s,v)$ for $v=12$ and $s=3$. \protect{\label{3,12}} }
\end{table} 

\begin{table}
$$
\begin{array}{|c||c|c|c|c|c|c|c|c||c|}\hline
\text{No.}& 0\dkrop 1& 1\dkrop 2&1\dkrop 3& 1\dkrop 4&2\dkrop 3&2\dkrop 4&3\dkrop 4&4\dkrop 5& e\\\hline
1 & 1 & 1 & 1 & 1 & 1 & 2 & 2 & 3 & 25 \\ \hline
2 & 1 & 1 & 1 & 1 & 1 & 1 & 3 & 3 & 26 \\ \hline
3 & 1 & 2 & 0 & 1 & 1 & 3 & 1 & 3 & 27 \\ \hline
4 & 1 & 2 & 1 & 1 & 1 & 2 & 2 & 2 & 24 \\ \hline
5 & 1 & 2 & 0 & 1 & 2 & 2 & 2 & 2 & 25 \\ \hline
6 & 1 & 2 & 0 & 1 & 1 & 2 & 2 & 3 & 26 \\ \hline
7 & 1 & 2 & 0 & 2 & 1 & 2 & 1 & 3 & 26 \\ \hline
8 & 2 & 2 & 1 & 1 & 1 & 1 & 2 & 2 & 24 \\ \hline
\end{array} 
$$
\caption{$\hat{L}(s,v)$ for $v=12$ and $s=4$. \protect{\label{4,12}}}
\end{table}
\begin{table}
{\small
$$
\begin{array}{|c||c|c|c|c|c|c|c|c|c|c|c|c||l|}\hline
\text{No.} & 0\dkrop 1& 1\dkrop 2&1\dkrop 3& 1\dkrop 4& 1\dkrop 5&2\dkrop 3&2\dkrop 4& 2\dkrop 5&3\dkrop 4&3\dkrop 5& 4\dkrop 5&  5\dkrop 6&e \\\hline
  1 & 1 & 1 & 1 & 0 & 1 & 2 & 1 & 0 & 1 & 2 & 1 & 1& 21\\ \hline
  2 & 1 & 1 & 1 & 0 & 1 & 2 & 1 & 0 & 2 & 1 & 1 & 1& 21\\ \hline
  3 & 1 & 1 & 1 & 0 & 0 & 3 & 0 & 0 & 1 & 3 & 1 & 1& 26\\ \hline
  4 & 1 & 1 & 1 & 0 & 0 & 3 & 0 & 0 & 2 & 2 & 1 & 1& 26\\ \hline
  5 & 1 & 1 & 1 & 0 & 1 & 1 & 2 & 0 & 1 & 1 & 2 & 1& 20\\ \hline
  6 & 1 & 1 & 1 & 0 & 1 & 1 & 1 & 1 & 1 & 1 & 1 & 2& 20\\ \hline
  7 & 1 & 1 & 1 & 0 & 1 & 1 & 0 & 1 & 2 & 1 & 1 & 2& 21\\ \hline
  8 & 1 & 1 & 1 & 0 & 1 & 1 & 1 & 0 & 1 & 2 & 1 & 2& 21\\ \hline
  9 & 1 & 1 & 1 & 0 & 1 & 1 & 0 & 1 & 1 & 2 & 1 & 2& 22'\\ \hline
10 & 1 & 1 & 1 & 0 & 0 & 2 & 0 & 1 & 2 & 1 & 1 & 2& 22\\ \hline
11 & 1 & 1 & 1 & 0 & 1 & 2 & 0 & 0 & 2 & 1 & 1 & 2& 22\\ \hline
12 & 1 & 1 & 1 & 0 & 0 & 2 & 0 & 1 & 1 & 2 & 1 & 2& 23\\ \hline
13 & 1 & 1 & 1 & 0 & 1 & 2 & 0 & 0 & 1 & 2 & 1 & 2& 23\\ \hline
14 & 1 & 1 & 1 & 0 & 0 & 2 & 0 & 0 & 2 & 2 & 1 & 2& 24\\ \hline
15 & 1 & 1 & 1 & 0 & 0 & 2 & 0 & 0 & 1 & 3 & 1 & 2& 25\\ \hline
16 & 1 & 1 & 1 & 0 & 0 & 1 & 0 & 2 & 1 & 1 & 1 & 3& 22\\ \hline
17 & 1 & 1 & 1 & 0 & 1 & 1 & 0 & 1 & 1 & 1 & 1 & 3& 22\\ \hline
18 & 1 & 1 & 1 & 0 & 0 & 1 & 0 & 1 & 1 & 2 & 1 & 3& 23\\ \hline
19 & 1 & 1 & 1 & 0 & 0 & 1 & 0 & 0 & 1 & 3 & 1 & 3& 24\\ \hline
20 & 1 & 1 & 1 & 0 & 0 & 1 & 1 & 1 & 1 & 1 & 2 & 2& 21\\ \hline
21 & 1 & 1 & 1 & 0 & 1 & 1 & 1 & 0 & 1 & 1 & 2 & 2& 21\\ \hline
22 & 1 & 1 & 1 & 0 & 0 & 1 & 0 & 1 & 2 & 1 & 2 & 2& 22\\ \hline
23 & 1 & 1 & 1 & 0 & 0 & 1 & 0 & 0 & 2 & 2 & 2 & 2& 23\\ \hline
24 & 1 & 1 & 1 & 0 & 0 & 2 & 0 & 0 & 1 & 2 & 2 & 2& 23\\ \hline
25 & 1 & 1 & 1 & 0 & 0 & 1 & 0 & 1 & 1 & 1 & 2 & 3& 22\\ \hline
26 & 1 & 1 & 1 & 0 & 0 & 1 & 0 & 0 & 1 & 2 & 2 & 3& 23\\ \hline
27 & 1 & 1 & 1 & 0 & 0 & 1 & 0 & 0 & 1 & 1 & 3 & 3& 22\\ \hline
28 & 1 & 1 & 2 & 0 & 0 & 2 & 0 & 0 & 2 & 2 & 1 & 1& 26\\ \hline
29 & 1 & 2 & 1 & 0 & 0 & 1 & 0 & 2 & 1 & 1 & 1 & 2& 22\\ \hline
30 & 1 & 2 & 1 & 0 & 0 & 1 & 1 & 1 & 1 & 1 & 1 & 2& 21\\ \hline
31 & 1 & 2 & 1 & 0 & 1 & 1 & 1 & 0 & 1 & 1 & 1 & 2& 20\\ \hline
32 & 1 & 2 & 1 & 0 & 0 & 1 & 0 & 1 & 1 & 2 & 1 & 2& 22\\ \hline
33 & 1 & 2 & 1 & 0 & 0 & 1 & 0 & 1 & 2 & 1 & 2 & 1& 21\\ \hline
34 & 1 & 2 & 1 & 0 & 0 & 2 & 0 & 0 & 1 & 2 & 1 & 2& 23\\ \hline
35 & 1 & 2 & 1 & 0 & 0 & 1 & 0 & 1 & 1 & 1 & 2 & 2& 21\\ \hline
36 & 1 & 2 & 1 & 0 & 0 & 2 & 0 & 0 & 1 & 1 & 2 & 2& 21\\ \hline
37 & 2 & 2 & 1 & 0 & 0 & 1 & 0 & 0 & 1 & 1 & 2 & 2& 20\\ \hline
\end{array} 
$$}
\caption{$\hat{L}(s,v)$ for $v=12$ and $s=5$. \protect{\label{5,12}}}
\end{table}
\begin{table}
\footnotesize{
$$
\begin{array}{|c||c|c|c|c|c|c|c|c|c|c|c|c|c|c|c|c|c||l|}\hline
\text{No.} & 0\dkrop 1& 1\dkrop 2&1\dkrop 3& 1\dkrop 4& 1\dkrop 5&1\dkrop 6&2\dkrop 3&2\dkrop 4& 2\dkrop 5&2\dkrop 6&3\dkrop 4&3\dkrop 5&3\dkrop 6& 4\dkrop 5&4\dkrop 6&  5\dkrop 6&6\dkrop 7&e\\\hline

 1  &  1 & 1 & 1 & 0 & 0 & 0 & 1 & 1 & 0 & 0 & 1 & 1 & 1 & 1 & 1 & 1 & 1 &  19'\\ \hline
 2  &  1 & 1 & 1 & 0 & 0 & 0 & 1 & 1 & 0 & 1 & 1 & 1 & 0 & 1 & 1 & 1 & 1 &  18\\ \hline
 3  &  1 & 1 & 1 & 0 & 0 & 1 & 1 & 1 & 0 & 0 & 1 & 1 & 0 & 1 & 1 & 1 & 1 &  18\\ \hline
 4  &  1 & 1 & 1 & 0 & 0 & 0 & 1 & 0 & 0 & 0 & 2 & 1 & 1 & 1 & 1 & 1 & 1 &  20'\\ \hline
 5  &  1 & 1 & 1 & 0 & 0 & 0 & 1 & 0 & 0 & 1 & 2 & 1 & 0 & 1 & 1 & 1 & 1 &  19\\ \hline
 6  &  1 & 1 & 1 & 0 & 0 & 0 & 1 & 2 & 0 & 0 & 1 & 0 & 0 & 1 & 2 & 1 & 1 &  20\\ \hline
 7  &  1 & 1 & 1 & 0 & 0 & 0 & 1 & 1 & 0 & 0 & 1 & 0 & 1 & 1 & 2 & 1 & 1 &  19\\ \hline
 8  &  1 & 1 & 1 & 0 & 0 & 0 & 2 & 0 & 0 & 0 & 1 & 1 & 1 & 1 & 1 & 1 & 1 &  20'\\ \hline
9   &  1 & 1 & 1 & 0 & 0 & 0 & 1 & 1 & 0 & 0 & 1 & 0 & 1 & 2 & 1 & 1 & 1 &  19\\ \hline
10 &  1 & 1 & 1 & 0 & 0 & 0 & 1 & 1 & 0 & 0 & 2 & 0 & 0 & 1 & 2 & 1 & 1 &  21\\ \hline
11 &  1 & 1 & 1 & 0 & 0 & 0 & 1 & 0 & 0 & 0 & 2 & 0 & 1 & 1 & 2 & 1 & 1 &  20\\ \hline
12 &  1 & 1 & 1 & 0 & 0 & 0 & 1 & 0 & 0 & 0 & 2 & 0 & 1 & 2 & 1 & 1 & 1 &  20\\ \hline
13 &  1 & 1 & 1 & 0 & 0 & 0 & 1 & 0 & 0 & 0 & 3 & 0 & 0 & 1 & 2 & 1 & 1 &  22\\ \hline
14 &  1 & 1 & 1 & 0 & 0 & 0 & 1 & 0 & 0 & 0 & 1 & 0 & 2 & 1 & 1 & 1 & 2 &  20\\ \hline
15 &  1 & 1 & 1 & 0 & 0 & 0 & 1 & 1 & 0 & 0 & 1 & 0 & 1 & 1 & 1 & 1 & 2 &  19\\ \hline
16 &  1 & 1 & 1 & 0 & 0 & 0 & 1 & 0 & 0 & 1 & 1 & 0 & 1 & 1 & 1 & 1 & 2 &  19\\ \hline
17 &  1 & 1 & 1 & 0 & 0 & 0 & 1 & 1 & 0 & 1 & 1 & 0 & 0 & 1 & 1 & 1 & 2 &  18\\ \hline
18 &  1 & 1 & 1 & 0 & 0 & 0 & 1 & 0 & 0 & 0 & 1 & 1 & 1 & 1 & 1 & 1 & 2 &  19'\\ \hline
19 &  1 & 1 & 1 & 0 & 0 & 0 & 1 & 0 & 0 & 1 & 1 & 1 & 0 & 1 & 1 & 1 & 2 &  18\\ \hline
20 &  1 & 1 & 1 & 0 & 0 & 1 & 1 & 1 & 0 & 0 & 1 & 0 & 0 & 1 & 1 & 1 & 2 &  18\\ \hline
21 &  1 & 1 & 1 & 0 & 0 & 0 & 1 & 0 & 0 & 0 & 2 & 0 & 1 & 1 & 1 & 1 & 2 &  20\\ \hline
22 &  1 & 1 & 1 & 0 & 0 & 0 & 1 & 0 & 0 & 1 & 2 & 0 & 0 & 1 & 1 & 1 & 2 &  19\\ \hline
23 &  1 & 1 & 1 & 0 & 0 & 0 & 1 & 0 & 0 & 0 & 1 & 0 & 1 & 1 & 2 & 1 & 2 &  19\\ \hline
24 &  1 & 1 & 1 & 0 & 0 & 0 & 1 & 1 & 0 & 0 & 1 & 0 & 0 & 1 & 2 & 1 & 2 &  19\\ \hline
25 &  1 & 1 & 1 & 0 & 0 & 0 & 1 & 0 & 0 & 0 & 1 & 0 & 1 & 2 & 1 & 2 & 1 &  18\\ \hline
26 &  1 & 1 & 1 & 0 & 0 & 0 & 1 & 0 & 0 & 0 & 2 & 0 & 0 & 1 & 2 & 1 & 2 &  20\\ \hline
27 &  1 & 1 & 1 & 0 & 0 & 0 & 1 & 0 & 0 & 0 & 1 & 0 & 1 & 1 & 1 & 2 & 2 &  19\\ \hline
28 &  1 & 1 & 1 & 0 & 0 & 0 & 1 & 1 & 0 & 0 & 1 & 0 & 0 & 1 & 1 & 2 & 2 &  18\\ \hline
29 &  1 & 1 & 1 & 0 & 0 & 0 & 1 & 0 & 0 & 0 & 2 & 0 & 0 & 1 & 1 & 2 & 2 &  19\\ \hline
\end{array}
$$
}
\caption{$\hat{L}(s,v)$ for $v=12$ and $s=6$. \protect{\label{6,12}}}
\end{table}

\begin{table}
\footnotesize{
$$
\begin{array}{|c||c|c|c|c|c|c|c|c|c|c|c|c|c|c|c|c|c||l|}\hline
\text{No.} & 0 \dkrop 1 & 1 \dkrop 2 & 1 \dkrop 3 &  2 \dkrop 3 & 2 \dkrop 4 & 2 \dkrop 5 & 2 \dkrop 7 & 3 \dkrop 4 & 3 \dkrop 5 & 3 \dkrop 7 & 4 \dkrop 5 & 4 \dkrop 6 & 4 \dkrop 7 & 5 \dkrop 6 & 5 \dkrop 7 & 6 \dkrop 7 & 7 \dkrop 8 & e\\\hline
  1 & 1 & 1 & 1 & 1 & 1 & 0 & 0 & 1 & 0 & 0 & 1 & 0 & 1 & 1 & 1 & 1 & 1 & 16\\ \hline
  2 & 1 & 1 & 1 & 1 & 0 & 0 & 0 & 1 & 1 & 0 & 1 & 0 & 1 & 1 & 1 & 1 & 1 & 17'\\ \hline
  3 & 1 & 1 & 1 & 1 & 0 & 0 & 1 & 1 & 1 & 0 & 1 & 0 & 0 & 1 & 1 & 1 & 1 & 17'\\ \hline
  4 & 1 & 1 & 1 & 1 & 0 & 0 & 0 & 1 & 2 & 0 & 1 & 0 & 0 & 1 & 1 & 1 & 1 & 19'\\ \hline
  5 & 1 & 1 & 1 & 1 & 0 & 0 & 0 & 1 & 0 & 1 & 1 & 0 & 1 & 1 & 1 & 1 & 1 & 16\\ \hline
  6 & 1 & 1 & 1 & 1 & 0 & 0 & 0 & 1 & 1 & 1 & 1 & 0 & 0 & 1 & 1 & 1 & 1 & 18'\\ \hline
  7 & 1 & 1 & 1 & 1 & 0 & 0 & 0 & 1 & 0 & 1 & 1 & 1 & 0 & 1 & 1 & 1 & 1 & 16\\ \hline
  8 & 1 & 1 & 1 & 1 & 0 & 1 & 0 & 1 & 0 & 1 & 1 & 0 & 0 & 1 & 1 & 1 & 1 & 17\\ \hline
  9 & 1 & 1 & 1 & 1 & 0 & 0 & 1 & 1 & 0 & 0 & 1 & 0 & 0 & 1 & 1 & 1 & 2 & 16\\ \hline
10 & 1 & 1 & 1 & 1 & 0 & 0 & 0 & 1 & 1 & 0 & 1 & 0 & 0 & 1 & 1 & 1 & 2 & 17'\\ \hline
11 & 1 & 1 & 1 & 1 & 0 & 0 & 0 & 1 & 0 & 1 & 1 & 0 & 0 & 1 & 1 & 1 & 2 & 17\\ \hline
12 & 1 & 1 & 1 & 1 & 0 & 0 & 0 & 1 & 0 & 0 & 1 & 0 & 1 & 1 & 1 & 1 & 2 & 16\\ \hline
13 & 1 & 1 & 1 & 1 & 0 & 0 & 0 & 2 & 0 & 0 & 1 & 0 & 1 & 1 & 1 & 1 & 1 & 17\\ \hline
14 & 1 & 1 & 1 & 1 & 0 & 0 & 0 & 1 & 1 & 0 & 2 & 0 & 0 & 1 & 1 & 1 & 1 & 18'\\ \hline
15 & 1 & 1 & 1 & 1 & 0 & 0 & 0 & 1 & 0 & 1 & 2 & 0 & 0 & 1 & 1 & 1 & 1 & 17\\ \hline
16 & 1 & 1 & 1 & 1 & 0 & 0 & 0 & 1 & 1 & 0 & 1 & 0 & 0 & 1 & 2 & 1 & 1 & 18'\\ \hline
17 & 1 & 1 & 1 & 1 & 0 & 1 & 0 & 1 & 0 & 0 & 1 & 0 & 0 & 1 & 2 & 1 & 1 & 17\\ \hline
18 & 1 & 1 & 1 & 1 & 0 & 0 & 0 & 1 & 0 & 0 & 2 & 0 & 0 & 1 & 2 & 1 & 1 & 17\\ \hline
19 & 1 & 1 & 1 & 1 & 0 & 0 & 0 & 1 & 0 & 0 & 1 & 0 & 0 & 1 & 2 & 1 & 2 & 16\\ \hline
20 & 1 & 1 & 1 & 1 & 0 & 0 & 0 & 2 & 0 & 0 & 2 & 0 & 0 & 1 & 1 & 1 & 1 & 18\\ \hline
21 & 1 & 1 & 1 & 1 & 0 & 0 & 0 & 1 & 0 & 0 & 1 & 0 & 0 & 1 & 1 & 2 & 2 & 16\\ \hline
\end{array}
$$
}
\caption{$\hat{L}(s,v)$ for $v=12$ and $s=7$. Columns $1\dkrop 4$, $1\dkrop 5$, $1\dkrop 6$, $1\dkrop 7$,  $2\dkrop 6$, $3\dkrop 6$ are omitted; they contain only $0$ entries. \protect{\label{7,12}}}
\end{table}

\begin{table}
{\small
$$
\begin{array}{|c||c|c|c|c|c|c|c|c|c|c|c|c|c|c|c||c|}\hline
\text{No.} & 0 \dkrop 1 & 1 \dkrop 2 & 1 \dkrop 3 & 2 \dkrop 3 & 3 \dkrop 4 & 3 \dkrop 6 & 3 \dkrop 8 & 4 \dkrop 5 & 4 \dkrop 6 & 5 \dkrop 6 & 5 \dkrop 8 & 6 \dkrop 7 & 6 \dkrop 8 & 7 \dkrop 8 & 8 \dkrop 9 & e\\\hline
1 & 1 & 1 & 1 & 1 & 1 & 0 & 0 & 1 & 0 & 1 & 1 & 1 & 1 & 1 & 1& 14\\ \hline
2 & 1 & 1 & 1 & 1 & 1 & 0 & 1 & 1 & 0 & 1 & 0 & 1 & 1 & 1 & 1& 15\\ \hline
3 & 1 & 1 & 1 & 1 & 1 & 0 & 0 & 1 & 1 & 1 & 0 & 1 & 1 & 1 & 1& 15\\ \hline
4 & 1 & 1 & 1 & 1 & 1 & 1 & 0 & 1 & 0 & 1 & 0 & 1 & 1 & 1 & 1& 16\\ \hline
5 & 1 & 1 & 1 & 1 & 1 & 0 & 0 & 1 & 0 & 2 & 0 & 1 & 1 & 1 & 1& 15\\ \hline
\end{array}
$$
}
\caption{$\hat{L}(s,v)$ for $v=12$ and $s=8$. Columns containing  only $0$ entries are omitted. \protect{\label{8,12}}}
\end{table}
\begin{table}
{\small
$$
\begin{array}{|c||c|c|c|c|c|c|c|c|c|c|c|c||c|}\hline
\text{No.} &0\dkrop 1 &  1\dkrop 2 &  1\dkrop 3 &  2\dkrop 3 &  3\dkrop 4 &  4\dkrop 5 &  5\dkrop 6 &  6\dkrop 7 &  7\dkrop 8 &  7\dkrop 9 &  8\dkrop 9 &  9\dkrop 10 &  e\\ \hline
1 & 1 &      1 &      1 &      1 &      1 &      1 &      1 &      1 &      1 &      1 &      1 &      1 &      13\\ \hline
\end{array}
$$}
\caption{$\hat{L}(s,v)$ for $v=12$ and $s=9$. Columns containing  only $0$ entries are omitted. \protect{\label{9,12}}}
\end{table}
\begin{table}
{\small
$$
\begin{array}{|c||c|c|c|c|c|c|c|c|c|c|c|c|c||c|}\hline
\text{No.} & 0\dkrop 1& 1\dkrop 2&1\dkrop 3& 1\dkrop 4& 1\dkrop 5&2\dkrop 3&2\dkrop 4& 2\dkrop 5&3\dkrop 4&3\dkrop 5& 4\dkrop 5&  5\dkrop 6&e\\\hline

1 &  1 & 1 & 1 & 0 & 1 & 1 & 0 & 1 & 1 & 2 & 1 & 3 &  27\\ \hline
2 &  1 & 1 & 1 & 0 & 0 & 2 & 0 & 0 & 2 & 2 & 2 & 2 &  28\\ \hline
3 &  1 & 2 & 1 & 0 & 1 & 1 & 1 & 1 & 1 & 1 & 1 & 2 &  24\\ \hline
4 &  1 & 2 & 1 & 0 & 0 & 2 & 0 & 1 & 1 & 2 & 1 & 2 &  27\\ \hline
5 &  1 & 2 & 1 & 0 & 0 & 1 & 1 & 1 & 1 & 1 & 2 & 2 &  25\\ \hline

\end{array} 
$$}
\caption{$\hat{L}(s,v)$ for $v=13$ and $s=5$. \protect{\label{5,13}}}
\end{table}

\begin{table}
\footnotesize{
$$
\begin{array}{|c||c|c|c|c|c|c|c|c|c|c|c|c|c|c|c|c|c||c|}\hline
\text{No.} & 0\dkrop 1& 1\dkrop 2&1\dkrop 3& 1\dkrop 4& 1\dkrop 5&1\dkrop 6&2\dkrop 3&2\dkrop 4& 2\dkrop 5&2\dkrop 6&3\dkrop 4&3\dkrop 5&3\dkrop 6& 4\dkrop 5&4\dkrop 6&  5\dkrop 6&6\dkrop 7&e\\\hline
1 &  1 & 1 & 1 & 1 & 0 & 0 & 1 & 1 & 0 & 0 & 1 & 1 & 1 & 1 & 1 & 1 & 1 &  22\\ \hline
2 &  1 & 1 & 1 & 0 & 0 & 0 & 2 & 1 & 0 & 0 & 1 & 1 & 1 & 1 & 1 & 1 & 1 &  23\\ \hline
3 &  1 & 1 & 1 & 1 & 0 & 0 & 1 & 1 & 0 & 0 & 1 & 0 & 1 & 2 & 1 & 1 & 1 &  23\\ \hline
4 &  1 & 1 & 1 & 0 & 0 & 0 & 2 & 0 & 0 & 0 & 2 & 1 & 1 & 1 & 1 & 1 & 1 &  25\\ \hline
5 &  1 & 1 & 1 & 0 & 0 & 1 & 1 & 1 & 0 & 0 & 1 & 0 & 1 & 1 & 1 & 1 & 2 &  22\\ \hline
6 &  1 & 1 & 1 & 0 & 0 & 0 & 1 & 1 & 0 & 0 & 1 & 0 & 1 & 1 & 2 & 1 & 2 &  23\\ \hline
7 &  1 & 1 & 1 & 0 & 0 & 0 & 1 & 0 & 0 & 0 & 2 & 0 & 1 & 1 & 2 & 1 & 2 &  24\\ \hline
8 &  1 & 1 & 1 & 0 & 0 & 0 & 1 & 0 & 0 & 0 & 1 & 1 & 1 & 1 & 1 & 2 & 2 &  23\\ \hline
\end{array} 
$$}
\caption{$\hat{L}(s,v)$ for $v=13$ and $s=6$. \protect{\label{6,13}}}
\end{table}

\begin{table}
\footnotesize{
$$
\begin{array}{|c||c|c|c|c|c|c|c|c|c|c|c|c|c|c|c|c|c||c|}\hline
\text{No.} & 0 \dkrop 1 & 1 \dkrop 2 & 1 \dkrop 3 & 2 \dkrop 3 & 2 \dkrop 4 & 2 \dkrop 5 & 2 \dkrop 7 & 3 \dkrop 4 & 3 \dkrop 5 & 3 \dkrop 7 & 4 \dkrop 5 & 4 \dkrop 6 & 4 \dkrop 7 & 5 \dkrop 6 & 5 \dkrop 7 & 6 \dkrop 7 & 7 \dkrop 8&e\\\hline

  1 &  1 & 1 & 1 & 1 & 1 & 0 & 0 & 1 & 1 & 0 & 1 & 0 & 1 & 1 & 1 & 1 & 1&  20\\ \hline
  2 &  1 & 1 & 1 & 1 & 0 & 0 & 0 & 1 & 1 & 1 & 1 & 1 & 0 & 1 & 1 & 1 & 1&  20\\ \hline
  3 &  1 & 1 & 1 & 1 & 0 & 0 & 1 & 1 & 1 & 0 & 1 & 0 & 0 & 1 & 1 & 1 & 2&  20\\ \hline
  4 &  1 & 1 & 1 & 1 & 0 & 0 & 0 & 1 & 1 & 1 & 1 & 0 & 0 & 1 & 1 & 1 & 2&  21\\ \hline
  5 &  1 & 1 & 1 & 1 & 0 & 0 & 0 & 1 & 1 & 0 & 1 & 0 & 1 & 1 & 1 & 1 & 2&  20\\ \hline
  6 &  1 & 1 & 1 & 1 & 0 & 0 & 0 & 1 & 0 & 1 & 1 & 0 & 1 & 1 & 1 & 1 & 2&  20\\ \hline
  7 &  1 & 1 & 1 & 1 & 0 & 0 & 0 & 2 & 1 & 0 & 1 & 0 & 1 & 1 & 1 & 1 & 1&  21\\ \hline
  8 &  1 & 1 & 1 & 1 & 0 & 0 & 0 & 1 & 2 & 0 & 1 & 0 & 0 & 1 & 2 & 1 & 1&  23\\ \hline
  9 &  1 & 1 & 1 & 1 & 0 & 1 & 0 & 1 & 1 & 0 & 1 & 0 & 0 & 1 & 2 & 1 & 1&  22\\ \hline
10 &  1 & 1 & 1 & 1 & 0 & 0 & 0 & 1 & 1 & 0 & 2 & 0 & 0 & 1 & 2 & 1 & 1&  22\\ \hline
11 &  1 & 1 & 1 & 1 & 0 & 0 & 0 & 1 & 1 & 0 & 1 & 0 & 0 & 1 & 2 & 1 & 2&  21\\ \hline
12 &  1 & 1 & 1 & 1 & 0 & 0 & 0 & 2 & 1 & 0 & 2 & 0 & 0 & 1 & 1 & 1 & 1&  22\\ \hline
13 &  1 & 1 & 1 & 1 & 0 & 0 & 0 & 1 & 1 & 0 & 1 & 0 & 0 & 1 & 1 & 2 & 2&  20\\ \hline
\end{array}
$$
}
\caption{$\hat{L}(s,v)$ for $v=13$ and $s=7$. Columns containing  only $0$ entries  are omitted. \protect{\label{7,13}}}
\end{table}
\begin{table}
{\small
$$
\begin{array}{|c||c|c|c|c|c|c|c|c|c|c|c|c|c|c|c|c||c|}\hline
\text{No.} & 0 \dkrop 1 & 1 \dkrop 2 & 1 \dkrop 3 & 2 \dkrop 3 & 3 \dkrop 4 & 3 \dkrop 5 & 3 \dkrop 6 & 3 \dkrop 8 & 4 \dkrop 5 & 4 \dkrop 6 & 5 \dkrop 6 & 5 \dkrop 8 & 6 \dkrop 7 & 6 \dkrop 8 & 7 \dkrop 8 & 8 \dkrop 9 & e\\\hline
1 &  1 &  1 &  1 &  1 &  1 &  1 &  0 &  0 &  0 &  1 &  1 &  1 &  1 &  1 &  1 &  1 &  18\\ \hline
2 &  1 &  1 &  1 &  1 &  0 &  1 &  0 &  0 &  1 &  1 &  1 &  1 &  1 &  1 &  1 &  1 &  18\\ \hline
3 &  1 &  1 &  1 &  1 &  0 &  1 &  1 &  0 &  0 &  1 &  1 &  1 &  1 &  1 &  1 &  1 &  18\\ \hline
4 &  1 &  1 &  1 &  1 &  0 &  1 &  0 &  1 &  0 &  1 &  1 &  1 &  1 &  1 &  1 &  1 &  19\\ \hline
5 &  1 &  1 &  1 &  1 &  0 &  2 &  0 &  0 &  0 &  1 &  1 &  1 &  1 &  1 &  1 &  1 &  19\\ \hline

\end{array}
$$
}
\caption{$\hat{L}(s,v)$ for $v=13$ and $s=8$. Columns containing  only $0$ entries  are omitted. \protect{\label{8,13}}}
\end{table}
\begin{table}
{\small
$$
\begin{array}{|c||c|c|c|c|c|c|c|c|c|c|c|c|c|c||c|}\hline
\text{No.} &0\dkrop 1 &  1\dkrop 2 &  1\dkrop 3 &  2\dkrop 3 &  3\dkrop 4 &  3\dkrop 7 & 4\dkrop 5 &  5\dkrop 6 &  5\dkrop 7 & 6\dkrop 7 &  7\dkrop 8 &  7\dkrop 9 &  8\dkrop 9 &  9\dkrop 10 &  e\\ \hline
1 & 1 & 1 & 1 & 1 & 1 & 0 & 1 & 1 & 1 & 1 & 1 & 1 & 1 & 1 & 16\\ \hline
2 & 1 & 1 & 1 & 1 & 1 & 1 & 1 & 1 & 0 & 1 & 1 & 1 & 1 & 1 & 17\\ \hline
\end{array}
$$
}
\caption{$\hat{L}(s,v)$ for $v=13$ and $s=9$. Columns containing  only $0$ entries  are omitted.\protect{\label{9,13}}}
\end{table}
\clearpage
\section{Acknowledgements and concluding remarks}
We extensively used solvers for mixed integer programming in the early stages of our investigations. Our first  (incompressible) 12-clique of 3-intervals has been found using A. Makhorin's GLPK (\textit{Gnu Linear Programming Kit}) for Windows \cite{Machorin} .  Other incompressible cliques have been produced with the aid of SCIP Optimization Suite \cite{Achter,  Koch}.  Subsequently, we have employed Gurobi Optimizer  5.1 to perform preliminary verification that no additional incompressible  cliques exist. Then we  have written Python scripts  to make the results mathematically sound.

We wish to thank Professor J. Zaks for sending us reprints of his papers.

\skok
\noindent
{\small
J.B.:  Wydzia{\l} Matematyki, Informatyki i Ekonometrii, Uniwersytet Zielonog\'orski, ul. Podg\'orna 50, \\65-246 Zielona G\'ora, Poland\\
{\tt J.Bojarski@wmie.uz.zgora.pl}}

\skok
\noindent{\small
A.K.:  Wydzia{\l} Matematyki, Informatyki i Ekonometrii, Uniwersytet Zielonog\'orski, ul. Podg\'orna 50,\\ 65-246 Zielona G\'ora, Poland\\
{\tt A.Kisielewicz@wmie.uz.zgora.pl}
}
\skok
\noindent{\small
K.P.:  Wydzia{\l} Matematyki, Informatyki i Ekonometrii, Uniwersytet Zielonog\'orski, ul. Podg\'orna 50, \\
65-246 Zielona G\'ora, Poland\\
and\\
Wydzia{\l} Matematyki, Informatyki i Architektury Krajobrazu, Katolicki Uniwersytet Lubelski, Poland\\
{\tt K.Przeslawski@wmie.uz.zgora.pl}}

\end{document}